\documentclass[12pt]{article}
\usepackage{mathrsfs}
\usepackage{amssymb}
\usepackage{amsmath}
\usepackage{amsbsy}
\usepackage{epsfig}
\usepackage{enumerate}
\usepackage{bm}
\newtheorem{lemma}{Lemma}[section]
\newtheorem{theorem}{Theorem}[section]
\newtheorem{proposition}{Proposition}[section]

\newtheorem{definition}{Deriniton}[section]
\newtheorem{remark}{Remark}[section]

\newbox\TempBox \newbox\TempBoxA

\def\ep{\textsf{E}} 
\def\Sbep{\widehat{\mathbb E}} 
\def\cSbep{\widehat{\mathcal E}} 
\def\extSbep{\Sbep^{\ast}} 
\def\Capc{\mathbb V} 
\def\cCapc{\mathcal V} 

\def\underwiggle 1{
\ifmmode\setbox\TempBox=\hbox{$ 1$}\else\setbox\TempBox=\hbox{
1}\fi \setbox\TempBoxA=\hbox to \wd\TempBox{\hss\char'176\hss}
\rlap{\copy\TempBox}\smash{\lower9pt\hbox{\copy\TempBoxA}} }
\renewcommand{\baselinestretch}{1.5}

\begin{document}

\thispagestyle{empty}

\begin{center}
 { \LARGE\bf Self-normalized moderate deviation and laws of the iterated logarithm under G-expectation$^{\ast}$}
\end{center}

\begin{center} {\sc
Li-Xin Zhang\footnote{Research supported by Grants from the National Natural Science Foundation of China (No. 11225104)  and the Fundamental Research Funds for the Central Universities.
}
}\\
{\sl \small Department of Mathematics, Zhejiang University, Hangzhou 310027} \\
(Email:stazlx@zju.edu.cn)\\
\today
\end{center}

\renewcommand{\abstractname}{~}
\begin{abstract}
\centerline{\bf Abstract} The sub-linear expectation or called G-expectation  is a nonlinear expectation having advantage of modeling non-additive probability  problems and the volatility uncertainty in finance. Let $\{X_n;n\ge 1\}$ be a sequence of independent random variables in a sub-linear expectation space $(\Omega, \mathscr{H}, \Sbep)$. Denote $S_n=\sum_{k=1}^n X_k$ and $V_n^2=\sum_{k=1}^n X_k^2$. In this paper, a  moderate deviation for self-normalized sums, that is,
the asymptotic capacity  of the event $\{S_n/V_n \ge x_n \}$ for $x_n=o(\sqrt{n})$, is
 found both for  identically distributed random variables and independent but not necessarily identically distributed random variables.
As an applications, the self-normalized laws of the
iterated logarithm are obtained.


\end{abstract}

\baselineskip 22pt

\renewcommand{\baselinestretch}{1.5}



\section{ Introduction}\label{sectIntro}
\setcounter{equation}{0}
Let $\{X_n; n\ge 1\}$ be a sequence of independent and identically distributed random variables on a probability space $(\Omega, \mathscr{F}, P)$. Set
$S_n=\sum_{j=1}^n X_j$, $ V_n^2=\sum_{j=1}^n X_j^2. $
The well-known classical
Hartman-Wintner law of the iterated logarithm (LIL) says that, if $EX_1=0$ and $EX_1^2=\sigma^2>0$, then
\begin{equation}\label{classicalLIL1} P\left(\limsup_{n\to \infty} \frac{S_n}{\sqrt{2\sigma^2n\log\log n}}=1\right)=1,
\end{equation}
and its  converse is obtained by Strassen (1966).  Griffin and Kuelbs (1989) obtained a self-normalized law of
the iterated logarithm under the following condition
\begin{equation}\label{classicalLIL2}
\lim_{x\to \infty} \frac{x^2P\big(|X_1|\ge x\big)}{EX_1^2I\{|X_1|\le x\}}=0.
\end{equation}
That is, if $EX_1=0$ and (\ref{classicalLIL2}) is satisfied, then
\begin{align}\label{classicalLIL3}
& P\left(\limsup_{n\to \infty} \frac{S_n}{V_n\sqrt{2\log\log n}}=1\right)= E\left[I \left\{\limsup_{n\to \infty} \frac{S_n}{V_n\sqrt{2\log\log n}}=1\right\}\right]=1.
\end{align}
On the other hand, Shao (1997)'s self-normalized moderate deviations gives us the asymptotic probability of $P(S_n\ge x_n V_n)$ as follows.
If $E X_1=0$ and (\ref{classicalLIL2}) is satisfied, then for any real sequence $\{x_n\}$ with $x_n\to \infty$ and $x_n=o(\sqrt{n})$,
\begin{equation}\label{classicalMD}
\lim_{n\to \infty} x_n^{-2}  \ln P\big( S_n\ge x_n V_n)=\lim_{n\to \infty} x_n^{-2}  \ln E\left[I\big\{ S_n\ge x_n V_n\big\}\right]=-\frac{1}{2}.
\end{equation}
The result is closely related to the Cram\'er (1938) large deviation. It is
known [cf. Petrov (1975)] that
$$ \lim_{n\to \infty} x_n^{-2}\ln P\left(\frac{S_n}{\sqrt{n}}\ge x_n\right)=-\frac{1}{2}.$$
holds for any sequence of $\{x_n\}$ with $x_n\to \infty$ and $x_n=o(\sqrt{n})$ 	 if and only if
$E X_1= 0$, $E X_1^2=1$ and $E\exp\{ t_0|X_1|\}<\infty$ for some $t_0 >0$. The   self-normalized limit theorems put a totally
new countenance upon classical limit theorems.

The  purpose of this paper is to study the self-normalized moderate deviation and self-normalized law of the iterated logarithm  for random variables in a sub-linear expectation space.
The sub-linear expectation or called G-expectation  is a nonlinear expectation advancing
the notions of  g-expectations, backward stochastic differential equations  and providing  a   flexible framework to model non-additive probability  problems and the volatility uncertainty in finance.
Peng (2006, 2008a, 2008b) introduced a general framework of the sub-linear expectation  of random variables   by relaxing the linear property of   the classical linear expectation to the sub-additivity
and positive homogeneity (cf.  Definition~\ref{def1.1} below), and introduced the notions of  G-normal random variable, G-Brownian motion,  independent and identically distributed random variables etc under  the sub-linear expectations. The construction of sub-linear expectations on   the space of continuous
paths and discrete time paths can also be found in   Yan et al (2012) and   Nutz and  Handel (2013).  For basic  properties of the sub-linear expectations, one can refer to Peng (2008b,2009,2010a,etc).
 Under the sub-linear expectation, the central limit theorem  was first established by Peng (2008b), large deviations and moderate deviations were derived by Gao and Xu (2011, 2012), and
 the Hartman-Winter laws of the iterated logarithm  were recently established by Chen and Hu (2014) for bounded random variables. Even for bounded random variables in a sub-linear expectation space, the self-normalized laws of the iterated logarithm can not follow from the Hartman-Winter laws of the iterated logarithm directly because $V_n^2/n$ does not converge to a constant.  We will  show that (\ref{classicalLIL3}) and (\ref{classicalMD}) are also true for random variables in a sub-linear expectation space when the expectation $E$ being replaced by the sub-linear expectation $\Sbep$.  The main difficult for proving the results under the sub-linear expectation  is that we can not use the additivity of the probability and the expectation  which is  essential in the proof of  classical results under the classical linear expectation. For example, the simple facts
 $-\Sbep[X_1^2I\{|X_1|\le x\}=\Sbep[-X_1^2I\{|X_1|\le x\}]$ and
 $$ \Sbep[X]=\int_0^{\infty} \Sbep\left[I\{X\ge x\}\right]dx+\int_{-\infty}^0\Sbep\left[I\{X\le -x\}\right]dx$$
 are not true now, and
 the conjugate method [cf. (4.9) of Petrov (1965)] is not available. These are the main keys to establish (\ref{classicalLIL2}) in Shao (1997).
   Our  main results on the self-normalized law of the iterated logarithm and self-normalized  moderate deviations for independent and identically distributed random variables   will be given in Section \ref{sectMain}.
In the next section, we state basic settings in a sub-linear expectation space including, capacity, independence, identical distribution  etc. One can skip this section if he/she is familiar with these  concepts.    The proofs are given in  Section \ref{sectproof}, where a Bernstein's type inequality for the maximum sum of independent random variables is also established for proving the law of the iterated logarithm. In Section \ref{SectGauss}, we consider the special case of $G$-normal random variables. We prove a finer self-normalized law of the iterated logarithm which shows the results are the same under the upper capacity and the lower capacity generated by the sub-linear expectation. In Section \ref{SectNonID}, we give similar results for independent but not necessarily identically distributed random variables.

\section{Basic Settings}\label{sectBasic}
\setcounter{equation}{0}

We use the framework and notations of Peng (2008b). Let  $(\Omega,\mathcal F)$
 be a given measurable space  and let $\mathscr{H}$ be a linear space of real functions
defined on $(\Omega,\mathcal F)$ such that if $X_1,\ldots, X_n \in \mathscr{H}$  then $\varphi(X_1,\ldots,X_n)\in \mathscr{H}$ for each
$\varphi\in C_b(\mathbb R^n)\bigcup  C_{l,Lip}(\mathbb R^n)$,  where $C_b(\mathbb R^n)$ denote  the
space  of all bounded  continuous functions and $C_{l,Lip}(\mathbb R^n)$ denotes the linear space of (local Lipschitz)
functions $\varphi$ satisfying
\begin{eqnarray*} & |\varphi(\bm x) - \varphi(\bm y)| \le  C(1 + |\bm x|^m + |\bm y|^m)|\bm x- \bm y|, \;\; \forall \bm x, \bm y \in \mathbb R^n,&\\
& \text {for some }  C > 0, m \in \mathbb  N \text{ depending on } \varphi. &
\end{eqnarray*}
$\mathscr{H}$ is considered as a space of ``random variables''. In this case we denote $X\in \mathscr{H}$.
Further, we let $C_{b,Lip}(\mathbb R^n)$  denote  the
space  of all bounded and Lipschitz functions   on $\mathbb R^n$.

\subsection{Sub-linear expectation and capacity}
\begin{definition}\label{def1.1} A {\bf sub-linear expectation} $\Sbep$ on $\mathscr{H}$  is a functional $\Sbep: \mathscr{H}\to \overline{\mathbb R}$ satisfying the following properties: for all $X, Y \in \mathscr H$, we have
\begin{description}
  \item[\rm (a)] {\bf Monotonicity}: If $X \ge  Y$ then $\Sbep [X]\ge \Sbep [Y]$;
\item[\rm (b)] {\bf Constant preserving}: $\Sbep [c] = c$;
\item[\rm (c)] {\bf Sub-additivity}: $\Sbep[X+Y]\le \Sbep [X] +\Sbep [Y ]$ whenever $\Sbep [X] +\Sbep [Y ]$ is not of the form $+\infty-\infty$ or $-\infty+\infty$;
\item[\rm (d)] {\bf Positive homogeneity}: $\Sbep [\lambda X] = \lambda \Sbep  [X]$, $\lambda\ge 0$.
 \end{description}
 Here $\overline{\mathbb R}=[-\infty, \infty]$. The triple $(\Omega, \mathscr{H}, \Sbep)$ is called a sub-linear expectation space. Give a sub-linear expectation $\Sbep $, let us denote the conjugate expectation $\cSbep$of $\Sbep$ by
$$ \cSbep[X]:=-\Sbep[-X], \;\; \forall X\in \mathscr{H}. $$
\end{definition}

Next, we introduce the capacities corresponding to the sub-linear expectations.
Let $\mathcal G\subset\mathcal F$. A function $V:\mathcal G\to [0,1]$ is called a capacity if
$$ V(\emptyset)=0, \;V(\Omega)=1 \; \text{ and } V(A)\le V(B)\;\; \forall\; A\subset B, \; A,B\in \mathcal G. $$
It is called to be sub-additive if $V(A\bigcup B)\le V(A)+V(B)$ for all $A,B\in \mathcal G$  with $A\bigcup B\in \mathcal G$.

 Let $(\Omega, \mathscr{H}, \Sbep)$ be a sub-linear space,
and  $\cSbep $  be  the conjugate expectation of $\Sbep$. It is natural to define the capacity of a set $A$ to be the sub-linear expectation of
the   indicator function $I_A$ of $A$. However, $I_A$ may be not in $\mathscr{H}$. So, we denote a pair $(\Capc,\cCapc)$ of capacities by
$$ \Capc(A):=\inf\{\Sbep[\xi]: I_A\le \xi, \xi\in\mathscr{H}\}, \;\; \cCapc(A):= 1-\Capc(A^c),\;\; \forall A\in \mathcal F, $$
where $A^c$  is the complement set of $A$.
Then $\Capc$  is sub-additive and
\begin{equation}\label{eq1.3} \begin{matrix}
&\Capc(A)=\Sbep[I_A], \;\; \cCapc(A)= \cSbep[I_A],\;\; \text{ if } I_A\in \mathscr H\\
&\Sbep[f]\le \Capc(A)\le \Sbep[g], \;\;\cSbep[f]\le \cCapc(A) \le \cSbep[g],\;\;
\text{ if } f\le I_A\le g, f,g \in \mathscr{H}.
\end{matrix}
\end{equation}

Further, we define an extension of $\Sbep^{\ast}$ of $\Sbep$ by
$$ \Sbep^{\ast}[X]=\inf\{\Sbep[Y]:X\le Y, \; Y\in \mathscr{H}\}, \;\; \forall X:\Omega\to \mathbb R, $$
where $\inf\emptyset=+\infty$. Then
$$ \begin{matrix}
&\Sbep^{\ast}[X]=\Sbep[X]\; \text{ if } X\in \mathscr H, \;\; \;\; \Capc(A)=\Sbep^{\ast}[I_A], \\
&\Sbep[f]\le \Sbep^{\ast}[X]\le \Sbep[g]\;\;
\text{ if } f\le X\le g, f,g \in \mathscr{H}.
\end{matrix}
$$

Also, we define the  Choquet integrals/expecations $(C_{\Capc},C_{\cCapc})$  by
$$ C_V[X]=\int_0^{\infty} V(X\ge t)dt +\int_{-\infty}^0\left[V(X\ge t)-1\right]dt $$
with $V$ being replaced by $\Capc$ and $\cCapc$ respectively.
It can be verified that (c.f. Zhang (2014)), if $\lim_{c\to \infty}\Sbep[(|X|-c)^+]=0$, then
\begin{equation} \label{eq1.5} \Sbep[|X|]\le C_{\Capc}(|X|).
\end{equation}

\begin{definition}\label{def3.1}
\begin{description}
\item{\rm (I)} A sub-linear expectation $\Sbep: \mathscr{H}\to \mathbb R$ is called to be  countably sub-additive if it satisfies
\begin{description}
  \item[\rm (e)] {\bf Countable sub-additivity}: $\Sbep[X]\le \sum_{n=1}^{\infty} \Sbep [X_n]$, whenever $X\le \sum_{n=1}^{\infty}X_n$,
  $X, X_n\in \mathscr{H}$ and
  $X\ge 0, X_n\ge 0$, $n=1,2,\ldots$;
 \end{description}

\item{\rm (II)}  A function $V:\mathcal F\to [0,1]$ is called to be  countably sub-additive if
$$ V\Big(\bigcup_{n=1}^{\infty} A_n\Big)\le \sum_{n=1}^{\infty}V(A_n) \;\; \forall A_n\in \mathcal F. $$

\item{\rm (III)}  A capacity $V:\mathcal F\to [0,1]$ is called a continuous capacity if it satisfies
\begin{description}
  \item[\rm (III1) ] {\bf Ccontinuity from below}: $V(A_n)\uparrow V(A)$ if $A_n\uparrow A$, where $A_n, A\in \mathcal F$;
  \item[\rm (III2) ] {\bf Continuity from above}: $V(A_n)\downarrow  V(A)$ if $A_n\downarrow A$, where $A_n, A\in \mathcal F$.
\end{description}
\end{description}
\end{definition}

Note that if $V$ is a   countably sub-additive capacity, then
\begin{align*}
0\le V\left (\bigcap_{n=1}^{\infty}\bigcup_{i=n}^{\infty}A_i\right)
\le V\left (\bigcup_{i=n}^{\infty}A_i\right)\le \sum_{i=n}^{\infty}V\left (A_i\right).
\end{align*}
So, we have the following Borel-Cantelli's Lemma.
\begin{lemma} ({\em Borel-Cantelli's Lemma}) Let $\{A_n, n\ge 1\}$ be a sequence of events in $\mathcal F$.
Suppose that $V$ is a countably sub-additive capacity.   If $\sum_{n=1}^{\infty}V\left (A_n\right)<\infty$, then
$$ V\left (A_n\;\; i.o.\right)=0, \;\; \text{ where } \{A_n\;\; i.o.\}=\bigcap_{n=1}^{\infty}\bigcup_{i=n}^{\infty}A_i. $$
\end{lemma}

\subsection{Independence and distribution}
\begin{definition} {\em(Peng (2006, 2008b))}

\begin{description}
  \item[ \rm (i)] ({\bf Identical distribution}) Let $\bm X_1$ and $\bm X_2$ be two $n$-dimensional random vectors defined
respectively in sub-linear expectation spaces $(\Omega_1, \mathscr{H}_1, \Sbep_1)$
  and $(\Omega_2, \mathscr{H}_2, \Sbep_2)$. They are called identically distributed, denoted by $\bm X_1\overset{d}= \bm X_2$  if
$$ \Sbep_1[\varphi(\bm X_1)]=\Sbep_2[\varphi(\bm X_2)], \;\; \forall \varphi\in C_{l,Lip}(\mathbb R^n), $$
whenever the sub-expectations are finite. A sequence $\{X_n;n\ge 1\}$ of random variables is said to be identically distributed if $X_i\overset{d}= X_1$ for each $i\ge 1$.
\item[\rm (ii)] ({\bf Independence})   In a sub-linear expectation space  $(\Omega, \mathscr{H}, \Sbep)$, a random vector $\bm Y =
(Y_1, \ldots, Y_n)$, $Y_i \in \mathscr{H}$ is said to be independent to another random vector $\bm X =
(X_1, \ldots, X_m)$ , $X_i \in \mathscr{H}$ under $\Sbep$  if for each test function $\varphi\in C_{l,Lip}(\mathbb R^m \times \mathbb R^n)$
we have
$$ \Sbep [\varphi(\bm X, \bm Y )] = \Sbep \big[\Sbep[\varphi(\bm x, \bm Y )]\big|_{\bm x=\bm X}\big],$$
whenever $\overline{\varphi}(\bm x):=\Sbep\left[|\varphi(\bm x, \bm Y )|\right]<\infty$ for all $\bm x$ and
 $\Sbep\left[|\overline{\varphi}(\bm X)|\right]<\infty$.
 \item[\rm (iii)] ({\bf IID random variables}) A sequence of random variables $\{X_n; n\ge 1\}$
 is said to be independent, if $X_{i+1}$ is independent to $(X_1,\ldots, X_i)$ for each $i\ge 1$, and it is  said to be  identically distributed, if
 $X_i\overset{d}=X_1$  for each $i\ge 1$.
 \end{description}
\end{definition}

\section{Main results}\label{sectMain}
\setcounter{equation}{0}
If $x\in \mathbb R$, $A\subset \mathbb R$, then the distance from $x$ to $A$ is defined as
$$d(x,A)=\inf_{y\in A}|x-y|.$$
If $\{x_n\}$ is  a real sequence, then $C(\{x_n\})$ denotes its cluster set, that is, $C(\{x_n\})=\{y: \liminf_{n\to \infty}|x_n-y|=0\}$. We write $\{x_n\}\twoheadrightarrow A$ if both $\lim_{n\to \infty} d(x_n,A)= 0 $ and $C(\{x_n\})=A$.

Let $\{X, X_n; n\ge 1\}$ be a sequence of independent and identically distributed random variables in the sub-linear expectation space  $(\Omega, \mathscr{H}, \Sbep)$ with $\Sbep[X]=\cSbep[X]=0$. Denote
$$ S_n=X_1+\cdots+X_n, \;\; V_n^2=X_1^2+\cdots+X_n^2. $$
 Define $l(x)=\Sbep[X^2\wedge x^2]$.

Our  main result is   the following self-normalized  law of iterated logarithm (LIL).
\begin{theorem}\label{thSelfLIL}  Suppose
\begin{description}
  \item[\rm(I)]  $\Capc(|X|\ge x)=o\big(x^{-2}l(x)\big)$ as $x\to \infty$;
  \item[\rm(II)] $\limsup_{x\to \infty} \frac{\Sbep[X^2\wedge x^2]}{\cSbep[X^2\wedge x^2]}<\infty$  (say, $<r^2<\infty$);
  \item[\rm (III)] $\Sbep[(|X|-c)^+]\to 0$ as $c\to \infty$.
\end{description}
Then
\begin{equation}\label{eqthLIL.1}
\cCapc\left(\limsup_{n\to \infty} \frac{|S_n|}{V_n \sqrt{2\log\log n}}\le 1\right)=1
\end{equation}
when $\Capc$ is countably sub-additive; and
\begin{equation}\label{eqthLIL.2}
\Capc\left(\left\{ \frac{S_n}{V_n \sqrt{2\log\log n}}\right\}\twoheadrightarrow [-1,1] \right)=1
\end{equation}
when $\Capc$ is continuous.
\end{theorem}

The proof of Theorem \ref{thSelfLIL} is based on the the following self-normalized moderated deviation.

\begin{theorem}\label{thMd} Suppose conditions (I)-(III) in Theorem  \ref{thSelfLIL} are satisfied and that
\begin{description}
  \item[\rm (IV)]  $x_n \to \infty$ and $x_n=o(\sqrt{n})$ as $n\to \infty$.
  \end{description}
  Then
  \begin{equation}\label{eqthMd}
  \lim_{n\to \infty} x_n^{-2} \ln \Capc\left(S_n\ge x_n V_n\right)=-\frac{1}{2}.
  \end{equation}
  Further, it also holds that
  $$\limsup_{n\to \infty} x_n^{-2} \ln \Capc\left(S_n\ge x_n V_n\right)\le -\frac{1}{2} $$
  if the condition $\Sbep[X]=\cSbep[X]=0$ is replaced by $\Sbep[X]\le 0$.
\end{theorem}

\begin{remark} If $\Sbep[X^2]<\infty$, $\cSbep[X^2]>0$ and $\Sbep[(X^2-c)^2]\to 0$ as $c\to \infty$, then conditions (I)-(III) are satisfied. Note
$$ \Sbep[X^2\wedge x^2]\le \Sbep[X^2\wedge(kx)^2]\le \Sbep[X^2\wedge x^2]+k^2x^2\Capc(|X|>x),\;\; k\ge 1. $$
The condition (I) implies that $l(x)$ is slowly varying as $x\to \infty$, i.e., for all $c>0$, $l(cx)/l(x)\to 1$ as $x\to \infty$. Further
\begin{align*}
& \frac{\extSbep[X^2I\{|X|\le x\}]}{l(x)}\to 1,\;\; \Sbep[|X|^r\wedge x^r]=o(x^{r-2}l(x)), \; r>2,\\
& C_{\Capc}\big(|X|^rI\{|X|\ge x\}\big)=\int_{x^r}^{\infty} \Capc(|X|^r\ge y)dy =o(x^{2-r} l(x)),\;\; 0<r<2.
\end{align*}
When Conditions (I) and (III) are satisfied,
$$ \Sbep[(|X|-x)^+]\le \extSbep[|X|I\{|X|\ge x\}] \le C_{\Capc}\big(|X|I\{|X|\ge x\}\big)=o(x^{-1} l(x)). $$

\end{remark}

 \section{Proofs }\label{sectproof}
 \setcounter{equation}{0}
 Let
 \begin{align*}
 & b=\inf\{x\ge 0: l(x)>0\}, \\
 & z_n=\inf\Big\{s:s\ge b+1, \frac{l(s)}{s^2}\le \frac{x_n^2}{n}\Big\}.
 \end{align*}
 Then $z_n\to \infty$, $nl(z_n)=x_n^2z_n^2$.
We follow the main idea of Shao (1997,1999) and Jing, Shao and Wang (2003). The main difference is that $-\Sbep[\cdot]$ and $\Sbep[-\cdot]$ may be different and  the conjugate method  is not available.  The proof of Theorem \ref{thMd} will be completed via four propositions.
\begin{proposition}\label{prop1} We have
\begin{equation}
\Capc\Big(S_n \ge x_n V_n, V_n^2\ge 9 n l(z_n)\Big )\le \exp\Big\{ -x_n^2+o(x_n)\Big\}.
\end{equation}
\end{proposition}
{\bf Proof}. Observe that $\Sbep[X_1^2\wedge z_n^2]=l(z_n)$, and
\begin{align*}
&\Capc\left( S_n\ge x_n V_n, V_n^2\ge 9nl(z_n)\right)\\
\le &\Capc\left( \sum_{i=1}^n \big((-z_n)\vee X_i\wedge z_n\big)\ge x_nV_n/2,  V_n^2\ge 9nl(z_n)\right) \\
  &+ \Capc\left( \sum_{i=1}^n    X_iI\{|X_i|> z_n\}\ge x_nV_n/2  \right)\\
  \le & \Capc\left( \sum_{i=1}^n \big((-z_n)\vee X_i\wedge z_n\big)\ge\frac{3}{2}  x_n ( nl(z_n))^{1/2} \right) \\
  &+ \Capc\left( \sum_{i=1}^n    I\{|X_i|> z_n\}\ge \frac{1}{4}x_n  \right)\\
  =:& J_1+J_2.
\end{align*}
Let $b_n=1/z_n$. As for $J_1$, by noting
$$ e^x\le 1+x+\frac{x^2}{2}+\frac{|x|^3}{6} e^{|x|}, $$
 we have
\begin{align*}
J_1
\le & \exp\left\{-\frac{3}{2}   \frac{x_n( nl(z_n))^{1/2}}{z_n}\right\} \Sbep \exp\left\{\frac{1}{z_n} \sum_{i=1}^n \big((-z_n)\vee X_i\wedge z_n\big)\right\}\\
\le & \exp\left\{-\frac{3}{2}   x_n^2\right\}\left[ \Sbep \exp\left\{\frac{1}{z_n}  \big((-z_n)\vee X_i\wedge z_n\big)\right\}\right]^n\\
\le & \exp\left\{-\frac{3}{2}   x_n^2\right\}\left[ 1+\frac{\Sbep\left[\big((-z_n)\vee X_i\wedge z_n\big)\right]}{z_n}+\frac{1}{2}
\frac{\Sbep[X_1^2\wedge z_n^2]}{z_n^2}+ \frac{e}{6} \frac{\Sbep[|X_1|^3\wedge z_n^3]}{z_n^3}  \right]^n\\
\le & \exp\left\{-\frac{3}{2}   x_n^2\right\}\left[ 1+\frac{\Sbep[X_1]+\Sbep\left[(|X_1|- z_n)^+\right]}{z_n}+\frac{1}{2}
\frac{l(z_n)}{z_n^2}+ \frac{e}{6} \frac{o\big(z_nl(z_n)\big)}{z_n^3}  \right]^n\\
\le & \exp\left\{-\frac{3}{2}   x_n^2\right\}\left[ 1+   \frac{1}{2}
\frac{l(z_n)}{z_n^2}+   \frac{o\big( l(z_n)\big)}{z_n^2}  \right]^n \\
\le & \exp\left\{-\frac{3}{2}   x_n^2\right\}\exp\left\{   \frac{1}{2}
\frac{nl(z_n)}{z_n^2}+   \frac{o\big(n l(z_n)\big)}{z_n^2}  \right\}
\le \exp\{-x_n^2+o(x_n^2)\}.
\end{align*}
As for $J_2$, we have
\begin{align}\label{eqproofprop1.1}
J_2\le & \exp\left\{-\frac{t}{4}x_n^2\right\}\extSbep\exp\left\{ t \sum_{i=1}^n    I\{|X_i|> z_n\}\right\}\nonumber\\
\le &  \exp\left\{-\frac{t}{4}x_n^2\right\}\extSbep \left[ \prod_{i=1}^n \left(1+(e^t-1)  I\{|X_i|> z_n\}\right)\right]\nonumber\\
\le &  \exp\left\{-\frac{t}{4}x_n^2\right\}\Sbep \left[ \prod_{i=1}^n \left(1+(e^t-1)  h(|X_i|/z_n)\right)\right]\nonumber\\
= &  \exp\left\{-\frac{t}{4}x_n^2\right\} \prod_{i=1}^n \left(1+(e^t-1) \Sbep \left[h(|X_i|/z_n)\right]\right)\nonumber\\
\le &  \exp\left\{-\frac{t}{4}x_n^2\right\}  \prod_{i=1}^n \left(1+(e^t-1)  \Capc(|X_i|> z_n/2)\right) \nonumber\\
\le &\exp\left\{-\frac{t}{4}x_n^2+(e^t-1)  n\Capc(|X_1|> z_n/2)\right\}\nonumber\\
\le &\exp\left\{-\frac{t}{4}x_n^2+(e^t-1)  n\frac{o\big(l(z_n)\big)}{z_n^2}\right\}\le \exp\left\{-2x_n^2+o(x_n^2)\right\}
\end{align}
for $t$ being chosen large enough, where $h(\cdot)$ is a Lipschitz function such that $I\{x>1\}\le h(x)\le I\{x>1/2\}$.   The proof is completed. $\Box$

Let $\lambda$ and $\theta>0$ be two real numbers. Define
\begin{equation}\label{eqfunctionf} f(s)=e^{\lambda s -\theta s^2}.
\end{equation}
Then
\begin{align*}
& f^{\prime}(s)=f(s) (\lambda-2\theta s), \;\;  f^{\prime}(0)=\lambda\\
& f^{\prime\prime}(s)=f(s)\big( (\lambda-2\theta s)^2-2\theta\big),\;\; f^{\prime\prime}(0)=\lambda^2-2\theta,\\
& f^{\prime\prime\prime}(s)=f(s)(\lambda-2\theta s)\big( (\lambda-2\theta s)^2-6\theta\big),\\
& f^{(4)}(s)=f(s)\big( (\lambda-2\theta s)^4-12\theta (\lambda-2\theta s)^2+12\theta^2\big).
\end{align*}
It is easily verified that
$$ |f(s)|\le e^{\frac{\lambda^2}{4\theta}}, \;\; |f^{\prime\prime}(s)|\le 4  e^{-\frac{3}{2}}\theta e^{\frac{\lambda^2}{4\theta}},
\;\; |f^{\prime\prime\prime}(s)|\le 3 \theta^{3/2} e^{\frac{\lambda^2}{4\theta}}.$$
It follows that
\begin{align*}
 f(s)= & 1+\lambda s+\big(\frac{\lambda^2}{2}-\theta \big)(s^2\wedge 1)+g(s) \;\;\text{ with } \\
  |g(s)|\le &\frac{1}{2} \theta^{3/2} e^{\frac{\lambda^2}{4\theta}}(|s|\wedge 1)^3+|\lambda| (|s|-1)^++  e^{\frac{\lambda^2}{4\theta}}I\{|s|>1\}.
  \end{align*}
For     non-negative $b$ and a random variable $\xi$, we have
\begin{align*}
& \Sbep\exp\left\{\lambda(b\xi)-\theta(b\xi)^2\right\}=\Sbep\left[1+\lambda b\xi+\big(\frac{\lambda^2}{2}-\theta\big)\big((b\xi)^2\wedge 1\big)+g(b\xi)\right]\\
\le &1+\Sbep[\lambda b\xi]+\Sbep\left[\big(\frac{\lambda^2}{2}-\theta\big)\big((b\xi)^2\wedge 1\big)+g(b\xi)\right]\\
\le &1+\Sbep\left[\big(\frac{\lambda^2}{2}-\theta\big)\big((b\xi)^2\wedge 1\big)+g(b\xi)\right]\\
\le & 1+ \big(\frac{\lambda^2}{2}-\theta\big)^+\Sbep\left[\big((b\xi)^2\wedge 1\big)\right]-\big(\frac{\lambda^2}{2}-\theta\big)^-\cSbep\left[\big((b\xi)^2\wedge 1\big)\right]
+\Sbep[|g(b\xi)|]
\end{align*}
if $\Sbep[\xi]\le 0$, and
\begin{align*}
& \Sbep\exp\left\{\lambda(b\xi)-\theta(b\xi)^2\right\}\\
\ge &1+\cSbep[\lambda b\xi]+\Sbep\left[\big(\frac{\lambda^2}{2}-\theta\big)\big((b\xi)^2\wedge 1\big)+g(b\xi)\right]\\
\ge &  1+ \big(\frac{\lambda^2}{2}-\theta\big)^+\Sbep\left[\big((b\xi)^2\wedge 1\big)\right]-\big(\frac{\lambda^2}{2}-\theta\big)^-\cSbep\left[\big((b\xi)^2\wedge 1\big)\right]
-\Sbep[|g(b\xi)|]
\end{align*}
if $\cSbep[\xi]\ge 0$.  Hence we obtain the following lemma.
\begin{lemma}\label{lemMomentG1} Suppose that $b$ is a positive number and   $\xi$ is a random variable. Then
\begin{align*}
 & \Sbep\exp\left\{\lambda(b\xi)-\theta(b\xi)^2\right\}\\
\le &  1+ \big(\frac{\lambda^2}{2}-\theta\big)^+\Sbep\left[\big((b\xi)^2\wedge 1\big)\right]-\big(\frac{\lambda^2}{2}-\theta\big)^-\cSbep\left[\big((b\xi)^2\wedge 1\big)\right]\\
& + O_{\lambda,\theta} \left(\Sbep[|b\xi|^3\wedge 1]+\Sbep[(|b\xi|-1)^+]+\Capc(|b\xi|>1)\right)
\end{align*}
if $\Sbep[\xi]\le 0$, and
\begin{align*}
 & \Sbep\exp\left\{\lambda(b\xi)-\theta(b\xi)^2\right\}\\
\ge &  1+ \big(\frac{\lambda^2}{2}-\theta\big)^+\Sbep\left[\big((b\xi)^2\wedge 1\big)\right]-\big(\frac{\lambda^2}{2}-\theta\big)^-\cSbep\left[\big((b\xi)^2\wedge 1\big)\right]\\
& + O_{\lambda,\theta} \left(\Sbep[|b\xi|^3\wedge 1]+\Sbep[(|b\xi|-1)^+]+\Capc(|b\xi|>1)\right)
\end{align*}
if $\cSbep[\xi]\ge 0$,
where $|O_{\lambda,\theta}|\le \frac{1}{2} \theta^{3/2} e^{\frac{\lambda^2}{4\theta}}+|\lambda|+e^{\frac{\lambda^2}{4\theta}}$.
\end{lemma}
Choose $b=b_n=z_n^{-1}$ if $\frac{\lambda^2}{2}-\theta\ge 0$ and $b=\underline{b}_n=\frac{1}{z_n}\frac{\Sbep[X_1^2\wedge z_n^2]}{\cSbep[X_1^2\wedge z_n^2]}$ if $\frac{\lambda^2}{2}-\theta< 0$. Then $z_n^{-2}\le \underline{b}_n\le C z_n^{-1}$, and so
\begin{align*}
 & \Sbep\exp\left\{\lambda(bX_1)-\theta(bX_1)^2\right\}\\
= &  1+ \big(\frac{\lambda^2}{2}-\theta\big)^+b_n^2\Sbep\left[\big(X_1^2\wedge b_n^{-2}\big)\right]-\big(\frac{\lambda^2}{2}-\theta\big)^-\underline{b}_n^2\cSbep\left[\big(X_1^2\wedge \underline{b}_n^{-2}\big)\right]\\
& + O_{\lambda,\theta} \left(\Sbep[|bX_1|^3\wedge 1]+\Sbep[(|bX_1|-1)^+]+\Capc(|bX_1|>1)\right)\\
=& 1+\big(\frac{\lambda^2}{2}-\theta\big)\frac{l(z_n)}{z_n^2}+\frac{o(l(z_n))}{z_n^2}O_{\lambda,\theta}\\
=&1+\big(\frac{\lambda^2}{2}-\theta\big) \frac{x_n^2}{n}+o\big(\frac{x_n^2}{n}\big)O_{\lambda,\theta}=\exp\left\{\big(\frac{\lambda^2}{2}-\theta\big) \frac{x_n^2}{n}+o\big(\frac{x_n^2}{n}\big)O_{\lambda,\theta}\right\}.
\end{align*}
Hence we obtain the following lemma.
\begin{lemma} \label{lemMomentG2} Let $b=b_n=\frac{1}{z_n}$ if $\frac{\lambda^2}{2}-\theta\ge 0$ and $b=\underline{b}_n=\frac{1}{z_n}\frac{\Sbep[X_1^2\wedge z_n^2]}{\cSbep[X_1^2\wedge z_n^2]}$ if $\frac{\lambda^2}{2}-\theta< 0$. Then
$$ \Sbep\exp\left\{\lambda(bS_n)-\theta(bV_n)^2\right\}=
\exp\left\{\big(\frac{\lambda^2}{2}-\theta\big) x_n^2+o(x_n^2)O_{\lambda,\theta}\right\}. $$
Further, it also holds that
$$ \Sbep\exp\left\{\lambda(bS_n)-\theta(bV_n)^2\right\}\le
\exp\left\{\big(\frac{\lambda^2}{2}-\theta\big) x_n^2+o(x_n^2)O_{\lambda,\theta}\right\}  $$
if the condition $\Sbep[X]=\cSbep[X]=0$ is replaced by $\Sbep[X]\le 0$.
\end{lemma}

\begin{proposition} \label{prop2} For $0<\delta<1$,
$$ \Capc\big(S_n\ge x_n V_n, \delta nl(z_n)\le V_n^2\le 9 nl(z_n)\big)\le \exp\{-x_n^2+o(x_n^2)\}. $$
\end{proposition}
{\bf Proof.} Let $1<\theta<2$, $b=b_n=z_n^{-1}$. Then by Lemma \ref{lemMomentG2},
\begin{align*}
&  \Capc\left(S_n\ge x_n V_n, \delta nl(z_n)\le V_n^2\le 9 nl(z_n)\right) \\
\le &  \Capc\left(2\frac{x_n}{V_n}S_n-\big(\frac{x_n}{V_n} V_n\big)^2\ge x_n^2, \frac{1}{3} b_n\le \frac{x_n}{V_n}\le \delta^{-1/2} b_n\right)\\
\le & \Capc\left(\sup_{3^{-1}\le a\le \delta^{-1/2}}\left( 2a b_n S_n-(ab_nV_n)^2\right)\ge x_n^2\right)\\
\le &\sum_{\big[\frac{\log 3^{-1}}{\log \theta}\big]\le j\le \frac{\log \delta^{-1/2}}{\log \theta}} \Capc\left(\sup_{\theta^j\le a\le \theta^{j+1}}\left( 2a b_n S_n-(ab_nV_n)^2\right)\ge x_n^2\right)\\
\le &\sum_{\big[\frac{\log 3^{-1}}{\log \theta}\big]\le j\le \frac{\log \delta^{-1/2}}{\log \theta}}
 \Capc\left(  2\theta^{j+1} b_n S_n-(\theta^jb_nV_n)^2\ge x_n^2 \right)\\
 \le & \sum_{\big[\frac{\log 3^{-1}}{\log \theta}\big]\le j\le \frac{\log \delta^{-1/2}}{\log \theta}}
\exp\{-\frac{x_n^2}{2}\} \Sbep\exp\left\{\theta^{j+1} b_n S_n-\frac{1}{2}\theta^{2j}(b_nV_n)^2\right\}\\
\le & \sum_{\big[\frac{\log 3^{-1}}{\log \theta}\big]\le j\le \frac{\log \delta^{-1/2}}{\log \theta}}
\exp\{-\frac{x_n^2}{2}\} \exp\left\{ \left(\frac{1}{2}\theta^{2j+2} - \frac{1}{2}\theta^{2j}\right)x_n^2+o(x_n^2)\right\}\\
\le & \left(\frac{ \log \delta^{-1/2}  +\log 3}{\log \theta}+1\right)\exp\left\{-\frac{x_n^2}{2}+\frac{1}{2}(\theta^2-1)\delta^{-1}x_n^2+o(x_n^2)\right\}.
\end{align*}
Let $\theta^2=1+x_n^{-1}$. Then
\begin{align*}
& \Capc\left(S_n\ge x_n V_n, \delta nl(z_n)\le V_n^2\le 9 nl(z_n)\right)\\
 \le & \left(  \log \delta^{-1/2}+ 2\right)2x_n\exp\left\{-\frac{x_n^2}{2}+\frac{1}{2}\delta^{-1}x_n+o(x_n^2)\right\}\\
 \le & \exp\left\{-\frac{x_n^2}{2}+ o(x_n^2)\right\}. \;\; \Box
 \end{align*}

Next, we consider the event $\{S_n\ge x_n V_n, V_n^2\le \delta nl(z_n)\}$. We need  Bernstein's type inequalities.
 \begin{lemma}\label{LemBineq} Suppose that $\{Y_n\}$ is a sequence of random variables on $(\Omega,\mathscr{H}, \Sbep)$ with $|Y_n|\le a$. Set $T_n=Y_1+\cdots+Y_n$, $B_n^2=\sum_{i=1}^n \Sbep\left[(Y_i-\Sbep[Y_i])^2\right]$. Then
 \begin{align}\label{eqBIneq1}
& \Capc\left(\sum_{i=1}^n(Y_i-\Sbep[Y_i])\ge x\right)\nonumber\\
\le & \exp\left\{-\frac{x^2}{2(B_n^2+2ax)}\right\} \le  \exp\left\{-\frac{x^2}{8\sum_{i=1}^n \Sbep[Y_i^2]+4ax }\right\},
 \end{align}
  \begin{align}\label{eqBIneq2}
& \Capc\left(\max_{i\le n}T_i \ge x+ \max_{k\le n}  \sum_{i=k+1}^n \Sbep[-Y_i] + \sum_{i=1}^n \Sbep[Y_i] \right)\nonumber\\
\le & 4\exp\left\{-\frac{x^2}{2(B_n^2+2ax)}\right\} \le  4\exp\left\{-\frac{x^2}{8\sum_{i=1}^n \Sbep[Y_i^2]+4ax }\right\}
 \end{align}
 and \begin{align}\label{eqBIneq3}
 \cCapc\left(\sum_{i=1}^n(Y_i-\cSbep[Y_i])\ge x\right)
\le & \exp\left\{-\frac{x^2}{2(\sum_{i=1}^n \Sbep\left[\big(Y_i-\cSbep[Y_i]\big)^2\right]+2ax)}\right\}\nonumber\\
 \le  &\exp\left\{-\frac{x^2}{8\sum_{i=1}^n \Sbep[Y_i^2]+4ax }\right\}.
 \end{align}
 \end{lemma}

 {\bf Proof.} Suppose $0<\lambda\cdot 2 a\le c<1$. Then
\begin{align*} \Sbep\left[e^{\lambda(Y_i-\Sbep[Y_i])}\right]
 \le & \Sbep\left[ 1+\lambda(Y_i-\Sbep[Y_i])+\frac{\lambda^2}{2}(Y_i-\Sbep[Y_i])^2\big(1+\sum_{k=3}^{\infty} (2\lambda a)^{k-2}\big)\right]\\
 \le &1+ \lambda\Sbep[(Y_i-\Sbep[Y_i])]+\frac{\lambda^2\Sbep\left[(Y_i-\Sbep[Y_i])^2\right]}{2(1-c)}\\
 \le & \exp\left\{\frac{\lambda^2\Sbep\left[(Y_i-\Sbep[Y_i])^2\right]}{2(1-c)}\right\}.
 \end{align*}
 So,
  \begin{align*}  \Sbep\left[e^{\lambda\sum_{i=1}^n(Y_i-\Sbep[Y_i])}\right]=\prod_{i=1}^n  \Sbep\left[e^{\lambda(Y_i-\Sbep[Y_i])}\right]
 \le  \exp\left\{\frac{\lambda^2B_n^2}{2(1-c)}\right\}.
 \end{align*}
 Letting $\lambda=\frac{x}{B_n^2+2ax}$ and   $c=\frac{2ax}{B_n^2+2ax}$ yields
 \begin{align*}
   \Capc & \left(\sum_{i=1}^n(Y_i-\Sbep[Y_i])\ge x\right)\le   e^{-\lambda x}\Sbep\left[e^{\lambda\sum_{i=1}^n(Y_i-\Sbep[Y_i])}\right] \\
& \le   e^{-\lambda x}\exp\left\{\frac{\lambda^2B_n^2}{2(1-c)}\right\}\le \exp\left\{-\frac{x^2}{2(B_n^2+2ax)}\right\}.
 \end{align*}

 As for (\ref{eqBIneq2}), we first show that
 \begin{equation}\label{eqMaxMG}
 \Sbep\exp\left\{\lambda \max_{k\le n} \sum_{i=1}^k (Y_i-\cSbep[Y_i])\right\}\le 4 \Sbep\exp\left\{ \lambda \sum_{i=1}^n (Y_i-\cSbep[Y_i])\right\}.
 \end{equation}
 Let
 $$ Q_0=0,\;\; Q_k=\exp\left\{ \frac{\lambda}{2} \sum_{i=1}^k (Y_i-\cSbep[Y_i])\right\},\; \; M_k=\max_{i\le k}Q_i. $$
 Then $M_k\ge Q_k\ge 0$ and
 \begin{align*}
 Q_nM_n=\sum_{k=0}^{n-1} Q_{k+1}(M_{k+1}-M_k)+\sum_{k=0}^{n-1}(Q_{k+1}-Q_k)M_k.
 \end{align*}
 Note that $Q_{k+1}=M_{k+1}$ when $M_{k+1}\ne M_k$. So
 \begin{align*}
  Q_{k+1}(M_{k+1}-M_k)= M_{k+1}(M_{k+1}-M_k)\ge \frac{1}{2} M_{k+1}^2-\frac{1}{2}M_k^2.
 \end{align*}
 It follows that
 \begin{align*}
 Q_nM_n\ge  &\frac{1}{2}M_n^2 +\sum_{k=0}^{n-1} Q_kM_k\Big(e^{\frac{\lambda}{2}(Y_{k+1}-\cSbep[Y_{k+1}])}-1\Big)\\
 \ge & \frac{1}{2}M_n^2 +\frac{\lambda}{2}\sum_{k=0}^{n-1} Q_kM_k\big(  Y_{k+1}-\cSbep[Y_{k+1}] \Big).
 \end{align*}
 Hence
 \begin{align*}
 \frac{1}{2}&\Sbep[M_n^2]\le \Sbep[Q_nM_n]+\frac{\lambda}{2}\sum_{k=0}^{n-1} \Sbep[Q_kM_k]\Sbep\left[\cSbep[Y_{k+1}]-  Y_{k+1} \right]\\
 & = \Sbep[Q_nM_n]\le \Sbep\left[\frac{1}{4}  M_n^2+ Q_n^2\right]\le \frac{1}{4}\Sbep[M_n^2]+\Sbep[Q_n^2].
 \end{align*}
 Hence $\Sbep[M_n^2]\le 4\Sbep[Q_n^2]$. The proof of (\ref{eqMaxMG}) is now completed.

 Now, by (\ref{eqMaxMG}),
 \begin{align*}
 \Sbep\exp& \left\{\lambda\max_{k\le n} T_k\right\}\le   \exp\left\{\lambda\max_{k\le n} \sum_{i=1}^k\cSbep[Y_i] \right\}\Sbep\exp\left\{\lambda\max_{k\le n} \sum_{i=1}^k(Y_i-\cSbep[Y_i])\right\}\\
 \le & 4\exp\left\{\lambda\max_{k\le n}  \sum_{i=1}^k\cSbep[Y_i] \right\}\Sbep\exp\left\{ \lambda \sum_{i=1}^n(Y_i-\cSbep[Y_i])\right\}\\
 = & 4\exp\left\{ \lambda\max_{k\le n}  \sum_{i=1}^k\cSbep[Y_i]- \sum_{i=1}^n\cSbep[Y_i] +\lambda \sum_{i=1}^n\Sbep[Y_i] \right\}\Sbep\exp\left\{ \lambda \sum_{i=1}^n(Y_i-\Sbep[Y_i])\right\}\\
 \le & 4 \exp\left\{ \lambda\max_{k\le n}  \sum_{i=k+1}^n\Sbep[-Y_i] +\lambda  \sum_{i=1}^n\Sbep[Y_i] \right\}\exp\left\{\frac{\lambda^2B_n^2}{2(1-c)}\right\},
 \end{align*}
 when $2a\lambda\le c<1$. The remainder proof is the same as that of (\ref{eqBIneq1}).

 The proof of  (\ref{eqBIneq3}) is similar to that of ($\ref{eqBIneq1}$), if the following facts are noted:
 $$ \cCapc\left(\sum_{i=1}^n(Y_i-\cSbep[Y_i])\ge x\right)\le   e^{-\lambda x}\cSbep\left[e^{\lambda\sum_{i=1}^n(Y_i-\cSbep[Y_i])}\right]
  = e^{-\lambda x}\prod_{i=1}^n\cSbep\left[e^{\lambda(Y_i-\cSbep[Y_i])}\right] $$
  and
\begin{align*} \cSbep\left[e^{\lambda(Y_i-\cSbep[Y_i])}\right]
 \le & \cSbep\left[ 1+\lambda(Y_i-\Sbep[Y_i])+\frac{\lambda^2(Y_i-\Sbep[Y_i])^2}{2(1-c)}\right]\\
 \le &1+ \lambda\cSbep[(Y_i-\cSbep[Y_i])]+\frac{\lambda^2\Sbep\left[(Y_i-\Sbep[Y_i])^2\right]}{2(1-c)}, \;\; 0<2 a\lambda<1,
 \end{align*}
where the last inequality is  due the fact that $\cSbep[X+Y]\le \cSbep[X]+\Sbep[Y]$.
 $\Box$

 \begin{proposition}\label{prop3} Suppose $0<\delta<r^{-2}$. For $n$ large enough,
$$ \Capc\left( V_n^2\le \delta nl(z_n)\right)\le \exp\{-2x_n^2\}. $$
\end{proposition}
{\bf Proof}. Let $\varepsilon=r^{-2}-\delta$. Applying the Bernstein inequality (\ref{eqBIneq1}) again yields
\begin{align}\label{eqproofprop3.1}
&\Capc\left( V_n^2\le \delta nl(z_n)\right)\le \Capc\left( \sum_{i=1}^n X_i^2\wedge z_n^2 \le \delta nl(z_n)\right)\nonumber\\
\le &\Capc\left( \sum_{i=1}^n\big(- X_i^2\wedge z_n^2+\cSbep[X_i^2\wedge z_n^2]\big) \ge n\cSbep[X_1^2\wedge z_n^2]-\delta nl(z_n)\right)\nonumber\\
\le &\Capc\left(\sum_{i=1}^n\big(- X_i^2\wedge z_n^2+\cSbep[X_i^2\wedge z_n^2]\big) \ge r^{-2} n\Sbep[X_1^2\wedge z_n^2]-\delta nl(z_n)\right)\nonumber\\
= &\Capc\left( \sum_{i=1}^n\big(- X_i^2\wedge z_n^2+\cSbep[X_i^2\wedge z_n^2]\big) \ge \varepsilon nl(z_n)\right)\nonumber\\
\le & \exp\left\{-\frac{\big( \varepsilon nl(z_n)\big)^2}{8 n \Sbep[X_1^4\wedge z_n^4]+4z_n\cdot  \varepsilon  nl(z_n)}\right\}\nonumber\\
\le & \exp\left\{-\frac{\big( \varepsilon nl(z_n)\big)^2}{8 n\cdot o(z_n^2l(z_n))+4z_n\cdot \varepsilon nl(z_n)}\right\}\nonumber\\
\le & \exp\left\{-\frac{  \varepsilon nl(z_n) }{   o( z_n^2)+4z_n}\right\}=\exp\left\{-\frac{x_n^2}{o(1)}\right\}
\le \exp\{-2x_n^2\}. \;\; \Box
\end{align}

Now, the upper bound of (\ref{eqthMd}) follows from Propositions \ref{prop1}-\ref{prop3} immediately. As for the lower bound, we let $\underline{b}_n=\frac{1}{z_n}$ and $\eta_n=2\underline{b}_n S_n-(\underline{b}_n V_n)^2$. Note
$$ x_nV_n=\inf_b\frac{x_n^2+b^2V_n^2}{2b}\le \frac{x_n^2+\underline{b}_n^2V_n^2}{2\underline{b}_n}. $$
We have
$$ \Capc(S_n\ge x_nV_n)\ge   \Capc\left( S_n\ge \frac{x_n^2+\underline{b}_n^2V_n^2}{2\underline{b}_n}\right)\ge \Capc( 2\underline{b}_n  S_n-(\underline{b}_nV_n)^2\ge x_n^2). $$
The  lower bound of (\ref{eqthMd}) follows from the following proposition.
\begin{proposition}\label{prop4} For any $0<\beta<1$,
\begin{equation}\label{eqprop4.1}
\liminf_{n\to \infty} x_n^{-2} \ln \Capc\left( \eta_n\ge (1+\beta)x_n^2\right)\ge -\frac{(1+\beta/4)^2}{2}.
\end{equation}
\end{proposition}

{\bf Proof.}  By Lemma \ref{lemMomentG2},
$$ \Sbep \left[e^{\lambda \eta_n}\right]=\exp\left\{(2\lambda^2-\lambda)x_n^2 +o(x_n^2) \right\}, \;\; \lambda>1/2. $$
That is,
$$ \varphi(\lambda)=\lim_{n\to\infty} x_n^{-2}\ln \Sbep \left[e^{\lambda \eta_n}\right]= 2\lambda^2-\lambda, \;\; \lambda>1/2. $$
On the other hand, note that $\Capc$ is sub-additive and $\Sbep\big[\big(e^{\lambda\eta_n}-c\big)^+\big]\le c^{-1} \Sbep[e^{2\lambda\eta_n}]\to 0$ as $c\to\infty$. We have
\begin{align}\label{eqproofprop4.1}
&\Sbep\left[e^{\lambda\eta_n}\right]\le  \int_0^{\infty} \Capc\left( e^{\lambda\eta_n}>t\right)dt \\
\le &1+ \int_1^{\infty}\frac{\Sbep\left[e^{(1+\epsilon)\lambda\eta_n}\right]}{t^{1+\epsilon}}dt
\le  1+  \frac{1}{\epsilon} \Sbep\left[e^{(1+\epsilon)\lambda\eta_n}\right]. \nonumber
\end{align}
It follows that
$$ \lim_{n\to \infty} x_n^{-2} \ln \int_0^{\infty} \Capc\left( e^{\lambda\eta_n}>t\right)dt =2\lambda^2-\lambda, \;\; \lambda>1/2. $$
Let $G_n(t)=1-\Capc(\eta_n>t)$ and $F_n(t)=\lim_{y\searrow t}G_n(t)$. It is easy to verify that $F_n(t)$ is a probability distribution function.
Let $\xi_n$ be a random variable on the probability space $(\Omega, \mathscr{F},P)$ with the distribution $F_n(t)$. Then
\begin{align*}
\int_0^{\infty} \Capc\left( e^{\lambda\eta_n}>t\right)dt=& \int_{-\infty}^{\infty} \lambda e^{\lambda t}\Capc(\eta_n>t)dt\\
=& \int_{-\infty}^{\infty} \lambda e^{\lambda t}P(\xi_n>t)dt=E_P\left[e^{\lambda\xi_n}\right].
\end{align*}
 Hence
 $$ \lim_{n\to \infty} x_n^{-2}\ln E_P\left[e^{\lambda\xi_n}\right]=\varphi(\lambda)=2\lambda^2-\lambda,\;\; \lambda>1/2. $$
 Note for $x\in \{y: \varphi^{\prime}(\lambda)=y,   \exists\; \lambda>1/2\}=(1, \infty)$,
 $$I(x)=\sup_{\lambda>1/2} \{\lambda x-\varphi(\lambda)\}=\frac{(x+1)^2}{8}. $$
 By the G\" atner-Ellis Theorem (cf.  Dembo  and Zeitouni  (1998)),
 \begin{align*}
  &\liminf_{n\to \infty} x_n^{-2}\ln \Capc\left(\eta_n\ge (1+\beta)x_n^2\right) \\
  \ge & \liminf_{n\to \infty} x_n^{-2}\ln \left\{1- F_n(1+\beta/2)x_n^2) \right\}\\
 \ge & \liminf_{n\to \infty} x_n^{-2}\ln P\left(\xi_n/x_n^2>1+\beta/2\right)\\
 \ge & -\inf_{x\in (1+\beta/2,\infty)} I(x)=-\frac{(1+\beta/4)^2}{2}.\;\; \Box
 \end{align*}

 The proof of Theorem \ref{eqthMd} is now completed.

 \bigskip
 Now, we begin the proof of the self-normalized law of the iterated logarithm.

{\bf Proof of Theorem \ref{thSelfLIL}}. We first show (\ref{eqthLIL.1}). That is
\begin{equation}\label{eqproofLIL.1}
\limsup_{n\to\infty} \frac{|S_n|}{V_n\sqrt{2\log\log n}}\le 1\;\; a.s. \; \Capc.
\end{equation}
Let $m_k=[e^{k/(\log\log k)^2}]$, $x_k=(2\log\log m_k)^{1/2}$. Then $x_k\sim (2\log k)^{1/2}$.  Observer that for $0<\epsilon<1/2$,
 \begin{align}\label{eqproofLIL.2}
 &\Capc\left(\max_{m_k\le n\le m_{k+1}}\frac{S_n}{V_n}\ge (1+7\epsilon) x_k\right)\nonumber\\
 \le &\Capc\left( \frac{S_{m_k}}{V_{m_k}}\ge (1+2\epsilon) x_k\right)+
 \Capc\left(\max_{m_k\le n\le m_{k+1}}\frac{S_n-S_{m_k}}{V_{m_k}}\ge 5\epsilon  x_k\right).
 \end{align}
 By Theorem \ref{eqthMd},
  $$ \Capc\left( \frac{S_{m_k}}{V_{m_k}}\ge (1+2\epsilon) x_k\right)\le \exp\left\{ -(1+2\epsilon)x_k^2/2\right\}\le k^{-1-\epsilon}$$
  for every sufficiently large $k$. We estimate the second term in the right-hand side of (\ref{eqproofLIL.2}) below. Let
  $\overline{z}_k$ be the number such that
  $$ \frac{l(\overline{z}_k)}{\overline{z}_k^2}=\frac{x_k^2}{m_{k+1}-m_k}, $$
  and denote $T_n=\sum_{i=m_k+1}^n (-\overline{z}_k^2)\vee X_i \wedge \overline{z}_k^2$. Then
  $$ \frac{\overline{z}_kx_k\sqrt{m_{k+1}l(\overline{z}_k)}}{x_k^2 \overline{z}_k^2}=\frac{\sqrt{m_{k+1}}}{\sqrt{m_{k+1}-m_k}}\to \infty. $$
  Let $0<\delta^2<r^{-2}/4$. Observer that
  \begin{align*}
   &\Capc\left(\max_{m_k\le n\le m_{k+1}}\frac{S_n-S_{m_k}}{V_{m_k}}\ge 5\epsilon  x_k\right)\\
   \le & \Capc\left(\max_{m_k\le n\le m_{k+1}}T_n\ge 2\epsilon  x_kV_{m_k},V_{m,k}^2> \delta^2 m_{k+1}l(\overline{z}_k)\right) \\
   & + \Capc\left( V_{m_k}^2\le \delta^2m_{k+1}l(\overline{z}_k)\right)+\Capc\left(\sum_{n=m_k+1}^{m_{k+1}}(|X_i|-\overline{z}_k)^+\ge 3\epsilon  x_kV_{m_k} \right)\\
   \le & \Capc\left(\max_{m_k\le n\le m_{k+1}}T_n\ge 2\epsilon   \delta x_k\sqrt{ m_{k+1}l(\overline{z}_k)}\right) \\
   & + \Capc\left( V_{m_k}^2\le \delta^2m_{k+1}l(\overline{z}_k)\right)+\Capc\left(\sum_{n=m_k+1}^{m_{k+1}}I\{|X_i|>\overline{z}_k\}\ge (\epsilon  x_k)^2 \right).
  \end{align*}
  Note that
  \begin{align*}
  \sum_{i=m_k+1}^{m_{k+1}} \Sbep\left[\big((-\overline{z}_k)\vee X_i\wedge \overline{z}_k\big)^2\right]
  =(m_{k+1}-m_k) l(\overline{z}_k)=x_k^2 \overline{z}_k^2,
  \end{align*}
   \begin{align*}
 & \sum_{i=m_k+1}^{m_{k+1}} \left|\cSbep\left[ (-\overline{z}_k)\vee X_i\wedge \overline{z}_k \right]\right|
  +\sum_{i=m_k+1}^{m_{k+1}} \left|\Sbep\left[ (-\overline{z}_k)\vee X_i\wedge \overline{z}_k \right]\right| \\
  \le & \sum_{i=m_k+1}^{m_{k+1}}\Sbep\left[(|X_1|-\overline{z}_k)^+\right]
  =(m_{k+1}-m_k) \frac{o\big(l(\overline{z}_k)}{\overline{z}_k}\\
  =&o(x_k^2 \overline{z}_k)=o\left(  x_k\sqrt{ m_{k+1}l(\overline{z}_k)}\right).
  \end{align*}
  So, by the Bernstein inequality (\ref{eqBIneq2}), for sufficiently large $k$,
  \begin{align*}
&  \Capc\left(\max_{m_k\le n\le m_{k+1}}T_n\ge 2\epsilon   \delta x_k\sqrt{ m_{k+1}l(\overline{z}_k)}\right)  \\
\le & 4 \exp\left\{-\frac{\big(\epsilon   \delta x_k\sqrt{ m_{k+1}l(\overline{z}_k)}\big)^2}{8x_k^2 \overline{z}_k^2+
4 \overline{z}_k\cdot\epsilon   \delta x_k\sqrt{ m_{k+1}l(\overline{z}_k)}}\right\}
 \\
\le & 4 \exp\left\{-\frac{ (\epsilon   \delta)^2 \frac{ m_{k+1}}{m_{k+1}-m_k}x_k^4 \overline{z}_k^2 }{8x_k^2 \overline{z}_k^2+
4 \cdot\epsilon   \delta x_k^2\overline{z}_k^2\sqrt{ \frac{ m_{k+1}}{m_{k+1}-m_k}} }\right\}\\
\le &\exp\{-2x_k^2\}\le k^{-2}.
  \end{align*}
Note $\delta^2<r^{-2}/4$.  Similar to (\ref{eqproofprop3.1}), applying the Bernstein inequality (\ref{eqBIneq1}) again yields
\begin{align*}\label{eqproofLIL.3}
&\Capc\left( V_{m_k}^2\le \delta^2m_{k+1}l(\overline{z}_k)\right)\le \Capc\left( \sum_{i=1}^{m_k} X_i^2\wedge \underline{z}_k^2 \le \delta^2m_{k+1}l(\overline{z}_k)\right)\nonumber\\
\le &\Capc\left( \sum_{i=1}^{m_k}\big(- X_i^2\wedge \underline{z}_k^2+\cSbep[X_i^2\wedge \underline{z}_k^2]\big) \ge m_k\cSbep[X_1^2\wedge \underline{z}_k^2]-\delta^2m_{k+1}l(\overline{z}_k)\right)\nonumber\\
\le &\Capc\left( \sum_{i=1}^{m_k}\big(- X_i^2\wedge \underline{z}_k^2+\cSbep[X_i^2\wedge \underline{z}_k^2]\big) \ge m_kr^{-2}\Sbep[X_1^2\wedge \underline{z}_k^2]-\delta^2m_{k+1}l(\overline{z}_k)\right)\nonumber\\
\le &\Capc\left( \sum_{i=1}^{m_k}\big(- X_i^2\wedge \underline{z}_k^2+\cSbep[X_i^2\wedge \underline{z}_k^2]\big) \ge  \delta^2m_{k+1}l(\overline{z}_k)\right)\nonumber\\
\le & \exp\left\{-\frac{\big( \delta^2m_{k+1}l(\overline{z}_k)\big)^2}{8 m_k \Sbep[X_1^4\wedge \overline{z}_k^4]+4\overline{z}_k\cdot  \delta^2  m_{k+1}l(\overline{z}_k)}\right\}\nonumber\\
\le & \exp\left\{-\frac{\big( \delta^2m_{k+1}l(\overline{z}_k)\big)^2}{8 m_k\cdot o( \overline{z}_k^2l(\overline{z}_k))+4\overline{z}_k\cdot \delta^2  m_{k+1}l(\overline{z}_k)}\right\}\nonumber\\
\le & \exp\left\{-\frac{  \delta^2 m_{k+1}l(\overline{z}_k) }{   o( \overline{z}_k^2)+4\overline{z}_k}\right\}
\le \exp\{-2x_k^2\}\le k^{-2}.
\end{align*}
  Finally, similar to  (\ref{eqproofprop1.1}) we have for sufficiently large $t$,
\begin{align*}
&\Capc\left(\sum_{n=m_k+1}^{m_{k+1}}I\{|X_i|>\overline{z}_k\}\ge (\epsilon  x_k)^2 \right)\\
\le & \exp\left\{-t(\epsilon x_k)^2 +(e^t-1)  \frac{(m_{k+1}-m_k)o(l(\overline{z}_k))}{\overline{z}_k^2}\right\} \\
\le & \exp\left\{ -t (\epsilon x_k)^2+ o(x_k^2)\right\}\le \exp\{-x_k^2\}\le k^{-2}.
\end{align*}
Combing the above inequalities yields
$$ \sum_k \Capc\left(\max_{m_k\le n\le m_{k+1}}\frac{S_n}{V_n}\ge (1+7\epsilon) x_k\right)\le \sum_k\big(k^{-1-\epsilon}+3 k^{-2}\big)<\infty.$$
Note the countable sub-additivity of $\Capc$. By the Borel-Cantelli lemma,
\begin{align*}
 &\Capc\left(\limsup_{n\to \infty} \frac{S_n}{V_n\sqrt{2\log\log n}}\ge 1+8\epsilon \right)\\
 \le &\Capc\left(\max_{m_k\le n\le m_{k+1}}\frac{S_n}{V_n}\ge (1+7\epsilon) x_k,\;\; i.o.\right)=0.
\end{align*}
Let $\{\epsilon_i\}$ be a sequence with $\epsilon_i\searrow 0$. Then
\begin{align*}
 &\Capc\left(\limsup_{n\to \infty} \frac{S_n}{V_n\sqrt{2\log\log n}}>1 \right)\\
 \le &\sum_{i=1}^{\infty}\Capc\left(\limsup_{n\to \infty} \frac{S_n}{V_n\sqrt{2\log\log n}}\ge 1+8\epsilon_i \right)=0
\end{align*}
by countable sub-additivity of $\Capc$. (\ref{eqproofLIL.1}) is proved.

As for (\ref{eqthLIL.2}), note
$$ \frac{|X_n|}{V_n\sqrt{2n\log\log n}}\le \frac{1}{\sqrt{2\log\log n}}\to 0. $$
By (\ref{eqthLIL.1}) and Proposition 2.1 of Griffin and Kuelbs (1989), it is sufficient to show that
\begin{equation}\label{eqproofLIL.5}
\Capc\left(\limsup_{n\to \infty} \frac{S_n}{V_n\sqrt{2\log\log n}}\ge 1\; \text{ and } \liminf_{n\to \infty} \frac{S_n}{V_n\sqrt{2\log\log n}}\le -1\right)=1.
\end{equation}
Let $n_k=[e^{k(\log k)^2}]$, $x_n=\sqrt{2\log\log n}$. Observer that $x_{n_k}^2\sim 2\log k$.
Next, we show that
\begin{equation}\label{eqproofLIL.6}
\lim_{k\to \infty} \frac{V_{n_{k}}}{V_{n_{k+1}}}=0\;\; a.s. \; \Capc.
\end{equation}
By Proposition \ref{prop3} we have for $\delta<r^{-2}$,
\begin{align*}
 \Capc\left( V_{n_k}^2\le \delta n_k l(z_{n_k})\right)
\le
  \exp\{-2x_{n_k}^2\}\le k^{-2}.
\end{align*}
By the Borel-Cantelli lemma, we have
\begin{equation}\label{eqproofLIL.7}\liminf_{k\to \infty} \frac{V_{n_k}^2}{n_k l(z_{n_k})}\ge \delta \; a.s. \; \Capc.
\end{equation}
Also, for sufficiently large $k$,
\begin{align*}
&\Capc\left(\sum_{i=1}^{n_k}X_i^2\wedge z_{n_k}^2\ge \epsilon n_{k+1}l(z_{n_{k+1}})\right)\\
\le & \Capc\left(\sum_{i=1}^{n_k}\big(X_i^2\wedge z_{n_k}^2-\Sbep[X_i^2\wedge z_{n_k}^2]\ge \epsilon n_{k+1}l(z_{n_{k+1}})-n_kl(z_{n_k})\right)\\
\le & \Capc\left(\sum_{i=1}^{n_k}\big(X_i^2\wedge z_{n_k}^2-\Sbep[X_i^2\wedge z_{n_k}^2]\ge \frac{\epsilon}{2} n_{k+1}l(z_{n_k}))\right)\\
\le & \exp\left\{ -\frac{\big(\frac{\epsilon}{2} n_{k+1}l(z_{n_k})\big)^2}{8 o(1) z_{n_k}^2 n_kl(z_{n_k})+z_{n_k}^2 \cdot \frac{\epsilon}{2} n_{k+1}l(z_{n_k})}\right\}\le \exp\{-x_{n_k}^2\}\le k^{-2}.
\end{align*}
Finally, let $a_k^2=\frac{n_{k+1}}{n_k}z_{n_k}^2$
\begin{align*}
&\Capc\left(\sum_{i=1}^{n_k}(X_i^2- z_{n_k}^2)^+\ge \epsilon n_{k+1}l(z_{n_{k+1}})\right)\\
\le &\Capc\left(\sum_{i=1}^{n_k}(X_i^2- z_{n_k}^2)^+\ge \epsilon n_{k+1}l(z_{n_k})\right)=\Capc\left(\sum_{i=1}^{n_k}(X_i^2- z_{n_k}^2)^+\ge \epsilon a_k^2x_{n_k}^2\right)\\
\le &\Capc\left(\sum_{i=1}^{n_k}(X_i^2- z_{n_k}^2)^+\ge \epsilon a_k^2x_{n_k}^2, \max_{i\le n_k}|X_i|\le a_k\right)
+n_k\Capc\left(|X_1|\ge a_k\right)\\
\le &\Capc\left(\sum_{i=1}^{n_k}I\{|X_i|> z_{n_k}\}\ge \epsilon x_{n_k}^2\right)
+n_k\Capc\left(|X_1|\ge a_k\right).
\end{align*}
Similar to  similar to  (\ref{eqproofprop1.1}) we have
$$\Capc\left(\sum_{i=1}^{n_k}I\{|X_i|> z_{n_k}\}\ge \epsilon x_{n_k}^2\right)\le \exp\{-x_{n_k}^2\}\le k^{-2}. $$
Observer that
\begin{align*}
n_k\Capc\left(|X_1|\ge a_k\right)=& o(1)n_k \frac{l(a_k)}{a_k^2}\le c \frac{n_k^2l(a_k)}{n_{k+1}z_{n_k}^2}
\le cx_{n_k}^2\frac{n_k}{n_{k+1}}\frac{l(a_k)}{l(z_{n_k})}\\
\le & c x_{n_k}^2\big(\frac{n_k}{n_{k+1}}\big)^{1-\epsilon}\le k^{-2}.
\end{align*}
Combing the above inequality and applying the Borel-Cantelli lemma yield
\begin{equation}\label{eqproofLIL.8}\lim_{k\to \infty} \frac{V_{n_k}^2}{n_{k+1} l(z_{n_{k+1}})}=0 \; a.s. \; \Capc.
\end{equation}
Now, (\ref{eqproofLIL.6}) follows from (\ref{eqproofLIL.7}) and (\ref{eqproofLIL.8}).

Hence, by (\ref{eqproofLIL.6}) and (\ref{eqproofLIL.1}) we have
\begin{align}\label{eqproofLIL.9}
& \limsup_{n\to \infty} \frac{S_n}{V_n\sqrt{2\log\log n}}\ge \limsup_{k\to \infty} \frac{S_{n_k}}{V_{n_k}\sqrt{2\log\log n_k}} \nonumber\\
\ge & \limsup_{k\to \infty} \frac{S_{n_k}-S_{n_{k-1}}}{V_{n_k}\sqrt{2\log\log n_k}} -\limsup_{k\to \infty} \frac{|S_{n_{k-1}|}}{V_{n_k}\sqrt{2\log\log n_k}}\nonumber\\
=&  \limsup_{k\to \infty} \Big(1-\frac{V_{n_{k-1}}^2}{V_{n_k}^2}\Big)^{1/2}\frac{S_{n_k}-S_{n_{k-1}}}{\big(V_{n_k}^2-V_{n_{k-1}}^2\big)^{1/2}\sqrt{2\log\log n_k}} \nonumber\\
 & -\limsup_{k\to \infty}\frac{V_{n_{k-1}}}{V_{n_k}} \frac{|S_{n_{k-1}|}}{V_{n_{k-1}}\sqrt{2\log\log n_k}}\\
 =&  \limsup_{k\to \infty} \frac{S_{n_k}-S_{n_{k-1}}}{\big(V_{n_k}^2-V_{n_{k-1}}^2\big)^{1/2}\sqrt{2\log\log n_k}} \;\;a.s. \; \Capc
\end{align}
and similarly,
\begin{align}\label{eqproofLIL.10}
 \limsup_{n\to \infty} \frac{-S_n}{V_n\sqrt{2\log\log n}}
\ge    \limsup_{k\to \infty} \frac{-(S_{n_k}-S_{n_{k-1}})}{\big(V_{n_k}^2-V_{n_{k-1}}^2\big)^{1/2}\sqrt{2\log\log n_k}}  \;\;a.s. \; \Capc.
\end{align}
Let $h(x)$ be a   non-increasing Lipschitz function such that $I\{x\ge 1-2\epsilon\}\ge h(x)\ge I\{x\ge 1-\epsilon\}$. Denote
\begin{align*}
\eta_{k,1}=& h\left(\frac{S_{n_k}-S_{n_{k-1}}}{\big(V_{n_k}^2-V_{n_{k-1}}^2\big)^{1/2}\sqrt{2\log\log n_k}}\right),\\
\eta_{k,2}=& h\left(-\frac{S_{n_k}-S_{n_{k-1}}}{\big(V_{n_k}^2-V_{n_{k-1}}^2\big)^{1/2}\sqrt{2\log\log n_k}}\right).
\end{align*}
It can be verified that $(\eta_{k,1},\eta_{k,2})$, $k=1,2,\ldots,$ are independent bounded random vectors under $\Sbep$. Let $x_k=\sqrt{2\log\log n_k}$. By Theorem \ref{thMd}, for sufficiently large $k$ we have
\begin{align*} \Sbep[\eta_{k,1}]\ge & \Capc\Big(\frac{S_{n_k}-S_{n_{k-1}}}{\big(V_{n_k}^2-V_{n_{k-1}}^2\big)^{1/2}\sqrt{2\log\log n_k}}
\ge   1-\epsilon\Big)\\
\ge &\exp\{-\frac{1}{2}(1-\epsilon)x_k^2\}\ge ck^{1-\epsilon/2}.
\end{align*}
It follows that
$$ \sum_{k=1}^{\infty} \Sbep[\eta_{k,1}]=\infty.$$
So, by the Bernstein inequality (\ref{eqBIneq3}) we have
\begin{align*}
&1- \Capc  \left(\sum_{k=1}^n \eta_{k,1}>\frac{1}{2}\sum_{k=1}^n \Sbep[\eta_{k,1}]\right)
= \cCapc  \left(\sum_{k=1}^n \eta_{k,1}\le \frac{1}{2}\sum_{k=1}^n \Sbep[\eta_{k,1}]\right)\\
& =\cCapc\left(\sum_{k=1}^n \big(\eta_{k,1}-\Sbep[\eta_{k,1}]\big)\ge \frac{1}{2}\sum_{k=1}^n \Sbep[\eta_{k,1}]\right)\\
& \le  \exp\left\{-\frac{ (\sum_{k=1}^n \Sbep[\eta_{k,1}])^2/4}{8\sum_{k=1}^n \Sbep[\eta_{k,1}^2]+4\cdot \sum_{k=1}^n \Sbep[\eta_{k,1}]/2}\right\}\\
&\le \exp\left\{-\frac{1}{40}\sum_{k=1}^n \Sbep[\eta_{k,1}]\right\}\to 0\;\; \text{ as } n\to \infty.
\end{align*}
By the continuity of $\Capc$,
\begin{align*}
\Capc  \left(\sum_{k=1}^{\infty} \eta_{k,1}=\infty\right)=&
\Capc  \left(\bigcup_{n=1}^{\infty}\bigcap_{n=N}^{\infty}\left\{\sum_{k=1}^n \eta_{k,1}>\frac{1}{2}\sum_{k=1}^n \Sbep[\eta_{k,1}]\right\}\right)\\
\ge & \limsup_{n\to \infty}\Capc  \left(\sum_{k=1}^n \eta_{k,1}>\frac{1}{2}\sum_{k=1}^n \Sbep[\eta_{k,1}]\right)=1.
\end{align*}
Similarly,
$$\Capc  \left(\sum_{k=1}^{\infty} \eta_{k,2}=\infty\right)=1. $$
Now, by the independence of $\{(\eta_{k,1},\eta_{k,2})\}$,
\begin{align*}
&\Capc  \left(\sum_{k=1}^{\infty} \eta_{k,1}=\infty \text{ and } \sum_{k=1}^{\infty} \eta_{k,2}=\infty \right)\\
=&  \lim_{N_1\to \infty}\lim_{M_1\to \infty}\Capc\left(\sum_{k=1}^{N_1} \eta_{k,1}\ge M_1 \text{ and } \sum_{k=N_1+1}^{\infty} \eta_{k,2}=\infty \right)\\
=& \lim_{N_2\to \infty} \lim_{M_2\to \infty}\lim_{N_1\to \infty}\lim_{M_1\to \infty}\Capc\left(\sum_{k=1}^{N_1} \eta_{k,1}\ge M_1 \text{ and } \sum_{k=N_1+1}^{N_2} \eta_{k,2}\ge M_2 \right)\\
\ge & \lim_{N_2\to \infty} \lim_{M_2\to \infty}\lim_{N_1\to \infty}\lim_{M_1\to \infty}
\Sbep\left[g\big( \sum_{k=1}^{N_1} \eta_{k,1}/M_1\big)g\big(\sum_{k=N_1+1}^{N_2} \eta_{k,2}/M_2\big)\right]\\
=&\lim_{N_2\to \infty} \lim_{M_2\to \infty}\lim_{N_1\to \infty}\lim_{M_1\to \infty}
\Sbep\left[g\big( \sum_{k=1}^{N_1} \eta_{k,1}/M_1\big)\right]\Sbep\left[g\big(\sum_{k=N_1+1}^{N_2} \eta_{k,2}/M_2\big)\right]\\
\ge & \lim_{N_2\to \infty} \lim_{M_2\to \infty}\lim_{N_1\to \infty}\lim_{M_1\to \infty}\Capc\left(\sum_{k=1}^{N_1} \eta_{k,1}\ge 2 M_1\right)\Capc\left(\sum_{k=N_1+1}^{N_2} \eta_{k,2}\ge 2M_2 \right)\\
=&\Capc  \left(\sum_{k=1}^{\infty} \eta_{k,1}=\infty \right)\Capc\left(\sum_{k=1}^{\infty} \eta_{k,2}=\infty \right)=1,
\end{align*}
where $g(x)$ is a Lipschitz function with $I\{x\ge 1\}\ge g(x)\ge I\{x\ge 2\}$. Combing the above inequality, (\ref{eqproofLIL.9}) and (\ref{eqproofLIL.10}) we obtain
\begin{align*}
& \Capc\left(\limsup_{n\to \infty} \frac{S_n}{V_n\sqrt{2\log\log n}}\ge 1-2\epsilon\; \text{ and } \liminf_{n\to \infty} \frac{S_n}{V_n\sqrt{2\log\log n}}\le -(1-2\epsilon)\right)\\
& \quad \ge   \Capc  \left(\sum_{k=1}^{\infty} \eta_{k,1}=\infty \text{ and } \sum_{k=1}^{\infty} \eta_{k,2}=\infty \right)= 1.
\end{align*}
By the continuity of $\Capc$, letting $\epsilon\to 0$ we obtain (\ref{eqproofLIL.5}).   $\Box$

\section{Normal random variables}\label{SectGauss}
\setcounter{equation}{0}
In this section, we show that in (\ref{eqthLIL.2}) $\Capc$ can be replaced by $\cCapc$ when the random variables are normal distributed, and so
$$\cCapc\left(\left\{ \frac{S_n}{V_n \sqrt{2\log\log n}}\right\}\twoheadrightarrow [-1,1] \right)=\Capc\left(\left\{ \frac{S_n}{V_n \sqrt{2\log\log n}}\right\}\twoheadrightarrow [-1,1] \right)=1. $$
Let $0<\underline{\sigma}\le \overline{\sigma}<\infty$ and $G(\alpha)=\frac{1}{2}(\overline{\sigma}^2 \alpha^+ - \underline{\sigma}^2 \alpha^-)$. $X$ is call a normal $N\big(0, [\underline{\sigma}^2, \overline{\sigma}^2]\big)$ distributed random variable  (write $X\sim N\big(0, [\underline{\sigma}^2, \overline{\sigma}^2]\big)$) under $\Sbep$, if for any bounded Lipschitz function $\varphi$, the function $u(x,t)=\Sbep\left[\varphi\left(x+\sqrt{t} X\right)\right]$ ($x\in \mathbb R, t\ge 0$) is the unique viscosity solution of  the following heat equation:
      $$ \partial_t u -G\left( \partial_{xx}^2 u\right) =0, \;\; u(0,x)=\varphi(x). $$
 \begin{theorem} \label{thLILNormal}  Let $\{X, X_n; n\ge 1\}$ be a sequence of independent and identically distributed normal random variables with
 $X_i\sim N\big(0, [\underline{\sigma}^2, \overline{\sigma}^2]\big)$) under $\Sbep$.  Suppose that $\Capc$ is continuous.
Then
\begin{equation}\label{eqthLILNormal}
\cCapc\left(\left\{ \frac{S_n}{V_n \sqrt{2\log\log n}}\right\}\twoheadrightarrow [-1,1] \right)=1.
\end{equation}
\end{theorem}

To prove Theorem \ref{thLILNormal}, we recall the definition of $G$-Brownian motion.
 Let $C[0,\infty)$ be a function space of continuous functions   on $[0,\infty)$ with the norm $\|x\|=\sum\limits_{k=1}^{\infty} 2^{-k}\sup\limits_{0\le t\le k}|x(t)|$ and $C_b\big(C[0,\infty)\big)$ is the set of bounded continuous  functions $h(x):C[0,\infty)\to \mathbb R$.
It is showed that  there is a sub-linear expectation space $\big(\widetilde{\Omega}, \widetilde{\mathscr{H}},\widetilde{\mathbb E}\big)$ with
$\widetilde{\Omega}= C[0,\infty)$ and $C_b\big(C[0,\infty)\big)\subset \widetilde{\mathscr{H}}$ such that $(\widetilde{\mathscr{H}}, \widetilde{\mathbb E}[\|\cdot\|])$ is a Banach space, and
the canonical process $W(t)(\omega) = \omega_t  (\omega\in \widetilde{\Omega})$ is a G-Brownian motion with
$W(1)\sim N\big(0, [\underline{\sigma}^2, \overline{\sigma}^2]\big)$ under $\widetilde{\mathbb E}$, i.e.,
  for all $0\le t_1<\ldots<t_n$, $\varphi\in C_{l,lip}(\mathbb R^n)$,
\begin{equation}\label{eqBM} \widetilde{\mathbb E}\left[\varphi\big(W(t_1),\ldots, W(t_{n-1}), W(t_n)-W(t_{n-1})\big)\right]
  =\widetilde{\mathbb E}\left[\psi\big(W(t_1),\ldots, W(t_{n-1})\big)\right],
  \end{equation}
  where $\psi(x_1,\ldots, x_{n-1})\big)=\widetilde{\mathbb E}\left[\varphi\big(x_1,\ldots, x_{n-1}, \sqrt{t_n-t_{n-1}}W(1)\big)\right]$
  (c.f. Peng (2006, 2008a, 2010), Denis, Hu and Peng (2011)).

We denote   a pair of   capacities corresponding to the sub-linear expectation $\widetilde{\mathbb E}$ by  $(\widetilde{\Capc},\widetilde{\cCapc})$. Then $\widetilde{\Capc}$ and $\widetilde{\cCapc}$ are continuous.

Denis, Hu and Peng (2011)   showed the following representation of the G-Brownian motion (c.f, Theorem 52).
 \begin{lemma} \label{DenisHuPeng}
 Let $(\Omega, \mathscr{F}, P) $ be a probability
measure space and $\{B(t)\}_{t\ge 0}$  is a $P$-Brownian motion. Then for all bounded continuous function $\varphi: C_b[0,\infty)\to \mathbb R$,
$$ \widetilde{\mathbb E}\left[\varphi\big(W(\cdot)\big)\right]=\sup_{\theta\in \Theta}\ep_P\left[\varphi\big(W_{\theta}(\cdot)\big)\right],\;\;
W_{\theta}(t) = \int_0^t\theta(s) dB(s), $$
where
\begin{eqnarray*}
&\Theta=\left\{ \theta:\theta(t) \text{ is } \mathscr{F}_t\text{-adapted process such that }  \underline{\sigma}\le \theta(t)\le \overline{\sigma}\right\},&\\
& \mathscr{F}_t=\sigma\{B(s):0\le s\le t\}\vee \mathscr{N}, \;\; \mathscr{N} \text{ is the collection of } P\text{-null subsets}. &
\end{eqnarray*}
\end{lemma}

{\bf Proof of Theorem \ref{thLILNormal}}. By (\ref{eqthLIL.1})
 and Proposition 2.1 of Griffin and Kuelbs (1989), it is sufficient to show that
$$
\cCapc\left(\limsup_{n\to \infty} \frac{S_n}{V_n\sqrt{2\log\log n}}\ge 1\; \text{ and } \liminf_{n\to \infty} \frac{S_n}{V_n\sqrt{2\log\log n}}\le -1\right)=1.
$$
Note the sub-additive of $\Capc$ and the continuity of $\Capc$ and $\cCapc$. It is sufficient to show that for all $\epsilon>0$,
\begin{equation}
\cCapc\left(\limsup_{n\to \infty} \frac{S_n}{\epsilon+V_n\sqrt{2\log\log n}}> 1-\epsilon \right)=1.
\end{equation}
Let $W(t)$ be a $G$-Brwonian motion on $\big(\widetilde{\Omega}, \widetilde{\mathscr{H}},\widetilde{\mathbb E}\big)$ with
$W(1)\sim N\big(0, [\underline{\sigma}^2, \overline{\sigma}^2]\big)$. Denote $\widetilde{V}_n^2=\sum_{k=1}^n \left(W(k)-W(k-1)\right)^2$. By the continuity of   $\cCapc$ and $\widetilde{\cCapc}$ again,
\begin{align*}
&\cCapc\left(\limsup_{n\to \infty} \frac{S_n}{\epsilon+V_n\sqrt{2\log\log n}}> 1-\epsilon \right)\\
=& \lim_{n\to\infty} \lim_{N\to \infty} \cCapc\left(\max_{n\le k\le N} \frac{S_k}{\epsilon+V_k\sqrt{2\log\log k}}> 1-\epsilon \right)\\
\ge & \lim_{n\to\infty} \lim_{N\to \infty} \widetilde{\cCapc}\left(\max_{n\le k\le N} \frac{W(k)}{\epsilon+\widetilde{V}_k\sqrt{2\log\log k}}> 1-2\epsilon \right)\\
& \Big( \text{by the fact} \; \left( X_1,\ldots, X_N\right)\overset{d}=\left(W(1)-W(0), \ldots, W(N)-W(N-1)\right) \Big) \\
=&\widetilde{\cCapc}\left(\limsup_{n\to \infty} \frac{W(n)}{\epsilon+\widetilde{V}_n\sqrt{2\log\log n}}> 1-2\epsilon \right) \\
= &\inf_{\theta\in \Theta} P  \left( \limsup_{n\to \infty} \frac{W_{\theta}(n)}{\epsilon+V_{\theta}(n)\sqrt{2\log\log n}}> 1-2\epsilon \right)
\end{align*}
by Lemma \ref{DenisHuPeng}, where $V_{\theta}^2(n)=\sum_{k=1}^n \left(W_{\theta}(k)-W_{\theta}(k-1)\right)^2$.  So, it is sufficient to show that
for each $\theta\in \Theta$,
\begin{equation} \limsup_{n\to \infty} \frac{W_{\theta}(n)}{V_{\theta}(n)\sqrt{2\log\log n}}= 1 \;\; a.s. \; P.
\end{equation}
Let $m_k=\left(W_{\theta}(k)-W_{\theta}(k-1)\right)^2-\int_{k-1}^k\theta^2(s) ds$. It is easily seen that $\{m_k,\mathscr{F}_k\}$ is a sequence of martingale differences with $E_P[m_k^2|\mathscr{F}_{k-1}]\le 4 \overline{\sigma}^4$. By the law of large numbers for martingales,
$$\frac{1}{n}\left(V_{\theta}^2(n)-\int_0^n \theta^2(s)ds\right)= \frac{1}{n}\sum_{k=1}^n m_k \to 0 \; a.s.\;P. $$
It is obvious that $n\underline{\sigma}^2\le \int_0^n \theta^2(s)ds \le n\overline{\sigma}^2$. It follows that
$$ \frac{V_{\theta}^2(n)}{\int_0^n \theta^2(s)ds }\to 1\;\; a.s. \; P. $$
On the other hand, note that $W_{\theta}(t)=\int_0^t\theta(s)d B(s)$ is a continuous martingale with quadratic variation process $\langle W_{\theta},W_{\theta}\rangle (t)=\int_0^t\theta^2(s)ds$. By the Dambis-Dubins-Schwarz theorem, there is a standard Brownian motion $B$ under $P$ such  that   $W_{\theta}(t)=B\left(\langle W_{\theta},W_{\theta} \rangle_t\right)$.    So, it is sufficient to show that
\begin{equation}\label{eqproofthLILNormal.6} \limsup_{n\to \infty} \frac{B\left(\langle W_{\theta},W_{\theta} \rangle_n\right)}{\sqrt{2\langle W_{\theta},W_{\theta} \rangle_n \log\log n}}= 1 \;\; a.s. \; P.
\end{equation}
Note $\langle W_{\theta},W_{\theta} \rangle_t\to \infty$ and is a continuous function of $t$. By the law of the iterated logarithm for Brownian motion,
$$ \limsup_{t\to \infty} \frac{B\left(\langle W_{\theta},W_{\theta} \rangle_t\right)}{\sqrt{2\langle W_{\theta},W_{\theta} \rangle_t \log\log \langle W_{\theta},W_{\theta} \rangle_t}}=\limsup_{t\to \infty} \frac{B(t)}{\sqrt{2 t \log\log t}}= 1 \;\; a.s. \; P,
$$
which implies (\ref{eqproofthLILNormal.6}) by noting that $\langle W_{\theta},W_{\theta} \rangle_t \thickapprox t$, $\max\limits_{n\le t\le n+1}\big|\langle W_{\theta},W_{\theta} \rangle_t-\langle W_{\theta},W_{\theta} \rangle_n\big|\le \overline{\sigma}^2$ and the path properties of a Brownian motion.
The proof is now completed.    $\Box$

\begin{remark} We conjuncture that for all random variables satisfying the conditions Theorem \ref{thSelfLIL}, (\ref{eqthLILNormal}) holds, and
 \begin{equation}\label{eqthMdad}
  \lim_{n\to \infty} x_n^{-2} \ln \cCapc\left(S_n\ge x_n V_n\right)=\lim_{n\to \infty} x_n^{-2} \ln \Capc\left(S_n\ge x_n V_n\right)=-\frac{1}{2},
  \end{equation}
  whenever $x_n \to \infty$ and $x_n=o(\sqrt{n})$ as $n\to \infty$.
  (\ref{eqthLILNormal}) and (\ref{eqthMdad}) are interesting because they show that the self-normalized law of the iterated logarithm and the self-normalized moderate deviation under the sub-linear expectation are the same as those under the classical linear expectation. The self-normalization eliminates the effect of the non-linearity.

  Following the lines of the proofs of (\ref{eqproofLIL.5}) and Proposition \ref{prop4}, it is sufficient to show that there is a $\underline{b}_n>0$ such that
\begin{equation}\label{eqremark5.1.1}\lim_{n\to \infty} x_n^{-2} \ln \int_0^{\infty} \cCapc\left( e^{\lambda\eta_n}>t\right)dt =2\lambda^2-\lambda, \;\; \lambda>1/2,
\end{equation}
  where $\eta_n=2\underline{b}_n S_n-(\underline{b}_n V_n)^2$. With similar arguments as showing Lemma \ref{lemMomentG2} we can show that for
  $\underline{b}_n=\frac{1}{z_n}\frac{\Sbep[X_1^2\wedge z_n^2]}{\cSbep[X_1^2\wedge z_n^2]}$,
  $$ \cSbep \left[e^{\lambda \eta_n}\right]=\exp\left\{(2\lambda^2-\lambda)x_n^2 +o(x_n^2) \right\}, \;\; \lambda>1/2. $$
  Unfortunately, we are not able to conclude (\ref{eqremark5.1.1}) from the above equality  because  (\ref{eqproofprop4.1}) is not true for $\cSbep$.
\end{remark}

\section{Non-identically distributed random variables}\label{SectNonID}
\setcounter{equation}{0}
In this section, we consider the independent  but not necessarily identically distributed random variables. We give the self-normalized moderate deviation and the self-normalized law of the iterated logarithm similar to those for classical random variables in a probability space, which were established by Jing, Shao and Wang (2003).  Suppose that $\{X_n;n\ge 1\}$ is a sequence of independent random variables on the sub-linear expectation space $(\Omega, \mathscr{H}, \Sbep)$ with $\Sbep[X_k^2]<\infty$, $k=1,2,\ldots$. Let $\overline{B}_n^2=\sum_{k=1}^n \Sbep[X_k^2]$, $\underline{B}_n^2=\sum_{k=1}^n \cSbep[X_k^2]$, and
$$ \Delta_{n,x}=\frac{1}{\overline{B}_n^2}\sum_{k=1}^n \Sbep\left[X_k^2\left(1\wedge \Big|\frac{x}{\overline{B}_n}X_k\Big|\right)\right]. $$
\begin{theorem}\label{th6.1} Suppose $\Sbep[X_k]\le 0$. Let $q_n=\overline{B}_n^2/\underline{B}_n^2$. Then for $x\ge 2$,
\begin{equation} \label{eqth6.1.1}
\Capc\left(S_n\ge x V_n\right)
\le \exp\left\{ -\frac{x^2}{2}+O(1) q_n^3 \left(\log x+x^2\Delta_{n,x} \right)\right\},
\end{equation}
where $|O(1)|\le C$ with $C$ does not depend  on $x$.

Further, suppose $\Sbep[X_k]=\cSbep[X_k]=0$, $\limsup\limits_{n\to\infty}\overline{B}_n^2/\underline{B}_n^2<\infty$, $x_n\to \infty$ and
 $$ x_n^2\max_{i\le n}\Sbep[X_i^2]=o\left( \overline{B}_n^2\right), \;\; \Delta_{n,x_n}\to 0. $$
 Then
\begin{equation} \label{eqth6.1.3} \lim_{n\to \infty} x_n^{-2}\ln \Capc\left(S_n\ge x_n V_n\right)=-\frac{1}{2}.
\end{equation}
\end{theorem}

\begin{theorem}\label{th6.2}
 Suppose $\Sbep[X_k]=\cSbep[X_k]=0$, $\overline{B}_n\to \infty$, $\limsup\limits_{n\to\infty}\overline{B}_n^2/\underline{B}_n^2<\infty$,
 $$  \max_{i\le n}\Sbep[X_i^2]=o\left( \overline{B}_n^2/\log\log \overline{B}_n\right)    $$
 and that
 \begin{equation}\label{eqth6.2.2} \forall \epsilon>0, \;\; \frac{1}{\overline{B}_n^2}\sum_{i=1}^n \Sbep\left[\left(X_i^2-\epsilon \overline{B}_n^2/\log\log \overline{B}_n\right)^+\right]\to 0.
 \end{equation}
 Then
\begin{equation} \label{eqth6.2.3}
\cCapc\left(\limsup_{n\to \infty} \frac{|S_n|}{V_n \sqrt{2\log\log \overline{B}_n}}\le 1\right)=1
\end{equation}
when $\Capc$ is countably sub-additive; and
\begin{equation} \label{eqth6.2.4}
\Capc\left(\left\{ \frac{S_n}{V_n \sqrt{2\log\log\overline{B}_n}}\right\}\twoheadrightarrow [-1,1] \right)=1
\end{equation}
when $\Capc$ is continuous.

\end{theorem}

It is easily seen that the condition (\ref{eqth6.2.2}) implies that $\Delta_{n,x_n}\to 0$ for $x_n=(1\pm \epsilon)\sqrt{2\log\log \overline{B}_n^2}$ and then (\ref{eqth6.1.3}).   After having (\ref{eqth6.1.3}),  Theorem \ref{th6.1} can be proved by  similar arguments as showing  Theorem \ref{thSelfLIL} combing with the arguments as in the proof of Theorem 4.1 of Jing, Shao and Wang (2003). We omitted the details here.
The proof of Theorem \ref{th6.1} will be completed via four propositions.
\begin{proposition} \label{prop6.1} Suppose $\Sbep[X_k]\le 0$. We have
\begin{equation}
\Capc\left(S_n\ge x V_n, V_n^2\ge 9 \overline{B}_n^2\right)
\le 2\exp\left\{ -x^2+O(1)x^2 \Delta_{n,x}\right\}.
\end{equation}
\end{proposition}
{\bf Proof}. $b=b_x=x/\overline{B}_n$, $\widehat{S}_n=\sum_{i=1}^n X_i\wedge(A_0/b)$ where $A_0$ is an absolute constant to be determined later. Obser that
\begin{align*}
&\Capc\left(S_n\ge x V_n, V_n^2\ge 9 \overline{B}_n^2\right)\\
\le & \Capc\left(\widehat{S}_n\ge x V_n/2, V_n^2\ge 9 \overline{B}_n^2\right)+\Capc\left(\sum_{i=1}^n \big(X_i-A_0/b\big)^+\ge xV_n/2\right)\\
\le &\Capc\left(\widehat{S}_n\ge \frac{3}{2} x \overline{B}_n\right)+\Capc\left(\sum_{i=1}^n I\{b X_i>A_0\}\ge \frac{x^2}{4}\right).
\end{align*}
Note $e^s\le 1+s+\frac{s^2}{2}+\frac{e^s}{6} (s^3\vee 0)$. We have
\begin{align*}
&\Sbep\left[\exp\left\{\frac{3}{2}(bX_i)\wedge A_0\right\}\right]\\
\le &\Sbep\left[   1+\frac{3}{2}(bX_i)\wedge A_0+ \frac{9}{8} \big((bX_i)\wedge A_0\big)^2+\frac{27 e^{3A_0/2}}{48}  |bX_i|^3\wedge A_0^3 \right]\\
\le &   1+\frac{3}{2}\Sbep[bX_i] + \frac{9}{8} \Sbep\left[(bX_i)^2\right]+\frac{27 e^{3A_0/2}}{48}A_0^3 \Sbep\left[ |bX_i|^3\wedge1 \right]\\
\le &\exp\left\{ \frac{9}{8} b^2\Sbep\left[X_i^2\right]+\frac{27 e^{3A_0/2}}{48}A_0^3 \Sbep\left[ |bX_i|^3\wedge1 \right]\right\}.
\end{align*}
It follows that
\begin{align*}
 &\Capc\left(\widehat{S}_n\ge \frac{3}{2} x \overline{B}_n\right)\le \exp\left\{-\frac{9}{4}x^2\right\}\Sbep\left[\exp\left\{\frac{3}{2} b\widehat{S}_n\right\}\right] \\
\le & \exp\left\{-\frac{9}{4}x^2\right\}\prod_{i=1}^n \Sbep\left[\exp\left\{\frac{3}{2}(bX_i)\wedge A_0\right\}\right]\\
\le &\exp\left\{-\frac{9}{4}x^2\right\}\exp\left\{\frac{9}{4}b^2\overline{B}_n^2+\frac{27 e^{3A_0/2}}{48}A_0^3 x^2 \Delta_{n,x}\right\}\\
=& \exp\left\{-\frac{9}{4}x^2 +\frac{27 e^{3A_0/2}}{48}A_0^3 x^2\Delta_{n,x}\right\}.
\end{align*}
On the other hand, let $h(x)$ be a Lipschitz function such that $I\{x>A_0\}\le h(x)\le I\{x>A_0/2\}$. Then
\begin{align*}
&\Capc\left(\sum_{i=1}^n I\{b X_i>A_0\}\ge \frac{x^2}{4}\right)
\le \Capc\left(\sum_{i=1}^n h(b X_i)\ge \frac{x^2}{4}\right) \\
\le & \exp\left\{-t \frac{x^2}{4}\right\}\prod_{i=1}^n \Sbep\left[\exp\left\{t\; h(bX_i)\right\}\right]\\
\le & \exp\left\{-t \frac{x^2}{4}\right\}\prod_{i=1}^n \extSbep\left[1+e^t I\{x>A_0/2\}\right]\\
\le & \exp\left\{-t \frac{x^2}{4}\right\}\prod_{i=1}^n\left(1+e^t \frac{4}{A_0^2} \Sbep[(bX_i)^2]\right)\\
\le & \exp\left\{-t \frac{x^2}{4} +e^t \frac{4}{A_0^2} b^2 \overline{B}_n^2 \right\} \\
\le &\exp\left\{-t \frac{x^2}{4} +  \frac{4 e^t x^2}{A_0^2}  \right\}\le \exp\{-x^2\}
\end{align*}
if we choose $t=5$ and $A_0=120$. The proof is completed. $\Box$

\bigskip
Let $\lambda>0$ and $\theta>0$ be two real numbers. Define
 $ f(s)=e^{\lambda s -\theta s^2}$ as in (\ref{eqfunctionf}).  Then
\begin{align*}
 f(s)= & 1+\lambda s+\big(\frac{\lambda^2}{2}-\theta \big)s^2+g(s) \;\;\text{ with } \\
  |g(s)|\le &\Big(\frac{1}{2} \theta^{3/2} e^{\frac{\lambda^2}{4\theta}}+2    e^{-\frac{3}{2}}\theta e^{\frac{\lambda^2}{4\theta}}\Big) (s^2\wedge |s|^3).
  \end{align*}
Similar to Lemma \ref{lemMomentG1}, we have the following lemma.
\begin{lemma} \label{lem6.1} Suppose that $b$ is a positive number and   $\xi$ is a random variable. Then
\begin{align*}
 & \Sbep\exp\left\{\lambda(b\xi)-\theta(b\xi)^2\right\}\\
\le &  1+ b^2\big(\frac{\lambda^2}{2}-\theta\big)^+\Sbep\left[ \xi^2\right]-b^2\big(\frac{\lambda^2}{2}-\theta\big)^-\cSbep\left[\xi^2\right]
  + O_{\lambda,\theta} \left(\Sbep[|b\xi|^3\wedge (b\xi)^2]\right)
\end{align*}
if $\Sbep[\xi]\le 0$, and
\begin{align*}
 & \Sbep\exp\left\{\lambda(b\xi)-\theta(b\xi)^2\right\}\\
\ge &  1+ b^2\big(\frac{\lambda^2}{2}-\theta\big)^+\Sbep\left[ \xi^2\right]-b^2\big(\frac{\lambda^2}{2}-\theta\big)^-\cSbep\left[\xi^2\right]
  + O_{\lambda,\theta} \left(\Sbep[|b\xi|^3\wedge (b\xi)^2]\right)
\end{align*}
if $\cSbep[\xi]\ge 0$,
where $|O_{\lambda,\theta}|\le \frac{1}{2} \theta^{3/2} e^{\frac{\lambda^2}{4\theta}}+2    e^{-\frac{3}{2}}\theta e^{\frac{\lambda^2}{4\theta}}$.
\end{lemma}

Hence we have the following lemma similar to Lemma \ref{lemMomentG2}.
\begin{lemma} \label{lem6.2}   Suppose $\Sbep[X_i^2]<\infty$, $i\ge 1$, $x\ge 2$.
\begin{description}
  \item[\rm (a)] Suppose $\Sbep[X_i]\le 0$ ($i\ge 1$). Let $b=b_n=x/\overline{B}_n$ if $\frac{\lambda^2}{2}-\theta\ge 0$. Then
$$ \Sbep\exp\left\{\lambda(bS_n)-\theta(bV_n)^2\right\}\le
\exp\left\{\big(\frac{\lambda^2}{2}-\theta\big) x^2+ O_{\lambda,\theta}x^2\Delta_{n,x}\right\}. $$
  \item[\rm (b)] Suppose $\Sbep[X_i]\le 0$ ($i\ge 1$).
Let $b=b_n=x/\underline{B}_n$ if $\frac{\lambda^2}{2}-\theta< 0$. Then
$$ \Sbep\exp\left\{\lambda(bS_n)-\theta(bV_n)^2\right\}\le
\exp\left\{\big(\frac{\lambda^2}{2}-\theta\big) x^2+ O_{\lambda,\theta} q_n^3 x^2\Delta_{n,x}\right\}. $$
  \item[\rm (c)]   Suppose $\Sbep[X_i]\ge 0$ ($i\ge 1$), $\overline{B}_n^2\to \infty$,  $\frac{\lambda^2}{2}-\theta> 0$, $x_n\ge 2$,
  $$x_n^2\max_{i\le n}\Sbep[X_i^2]=o\left( \overline{B}_n^2 \right). $$
Let $b=b_n=x_n/\overline{B}_n$.  Then
$$ \Sbep\exp\left\{\lambda(bS_n)-\theta(bV_n)^2\right\}\ge
\exp\left\{\big(\frac{\lambda^2}{2}-\theta\big) x_n^2+ O_{\lambda,\theta}  x_n^2\Delta_{n,x_n}\right\}. $$
\end{description}
Here $|O_{\lambda,\theta}|\le C( \theta^{3/2} e^{\frac{\lambda^2}{4\theta}}+ \theta e^{\frac{\lambda^2}{4\theta}})$.
\end{lemma}

 \begin{proposition} \label{prop6.2}  Suppose $\Sbep[X_i]\le 0$, $i\ge 1$, and $0<\delta\le \frac{1}{4}\frac{\underline{B}_n^2}{\overline{B}_n^2}\le \frac{1}{4}$. For $x\ge 2$,
$$ \Capc\left(S_n\ge x V_n, \; V_n^2\le \delta \overline{B}_n^2\right)\le \exp\left\{-2x^2+O(1) x^2q_n^3 \Delta_{n,x}\right\}. $$
\end{proposition}
{\bf Proof}. Let $b=b_n=x/\underline{B}_n$. By Lemma \ref{lem6.2} (b) we have for $\lambda=2$,
\begin{align*}
& \Capc\left(S_n\ge x V_n, \; V_n^2\le \delta \overline{B}_n^2\right)\\
= & \Capc\left(b S_n\ge x \sqrt{(b V_n)^2}, \; (bV_n)^2\le \delta \frac{\overline{B}_n^2}{\underline{B}_n^2} x^2\right)\\
 \le &\Capc\left(b S_n\ge x \sqrt{(b V_n)^2}, \; (bV_n)^2\le \frac{1}{4} x^2\right)\\
 \le &\Capc\left(b S_n-2(bV_n)^2\ge 0 \right)
\le   \Sbep\exp\left\{\lambda\big(b S_n-2(bV_n)^2\big)\right\}\\
   \le &
\exp\left\{\big(\frac{\lambda^2}{2}-2\lambda\big) x^2+ O_{\lambda,2\lambda} q_n^3 x^2\Delta_{n,x}\right\}\\
= & \exp\left\{ -2 x^2+ O(1) q_n^3 x^2\Delta_{n,x}\right\}.\;\;\; \Box
\end{align*}

 \begin{proposition} \label{prop6.3}  Suppose $\Sbep[X_i]\le 0$, $i\ge 1$, and $ \delta= \frac{1}{4}\frac{\underline{B}_n^2}{\overline{B}_n^2}$. For $x\ge 2$,
$$ \Capc\left(S_n\ge x V_n, \; \delta \overline{B}_n^2\le V_n^2 \le 9 \overline{B}_n^2\right)\le \exp\left\{-\frac{x^2}{2}+O(1)q_n^{3/2} \left( \log x+ x^2 \Delta_{n,x}\right)\right\}. $$
\end{proposition}
{\bf Proof.} Let $1<\theta<2$, $b=b_n=\frac{x}{\overline{B}_n}$. Similar to the proof of Proposition \ref{prop2}, by Lemma \ref{lem6.2} (a) we have
\begin{align*}
&  \Capc\left(S_n\ge x V_n, \delta \overline{B}_n^2\le V_n^2\le 9 \overline{B}_n^2\right) \\
\le &  \Capc\left(2\frac{x}{V_n}S_n-\big(\frac{x}{V_n} V_n\big)^2\ge x^2, \frac{1}{3} b_n\le \frac{x}{V_n}\le \delta^{-1/2} b_n\right)\\
\le & \Capc\left(\sup_{3^{-1}\le a\le \delta^{-1/2}}\left( 2a b_n S_n-(ab_nV_n)^2\right)\ge x^2\right)\\
\le &\sum_{\big[\frac{\log 3^{-1}}{\log \theta}\big]\le j\le \frac{\log \delta^{-1/2}}{\log \theta}} \Capc\left(\sup_{\theta^j\le a\le \theta^{j+1}}\left( 2a b_n S_n-(ab_nV_n)^2\right)\ge x^2\right)\\
\le &\sum_{\big[\frac{\log 3^{-1}}{\log \theta}\big]\le j\le \frac{\log \delta^{-1/2}}{\log \theta}}
 \Capc\left(  2\theta^{j+1} b_n S_n-(\theta^jb_nV_n)^2\ge x^2 \right)\\
 \le & \sum_{\big[\frac{\log 3^{-1}}{\log \theta}\big]\le j\le \frac{\log \delta^{-1/2}}{\log \theta}}
\exp\{-\frac{x^2}{2}\} \Sbep\exp\left\{\theta^{j+1} b_n S_n-\frac{1}{2}\theta^{2j}(b_nV_n)^2\right\}\\
\le & \sum_{\big[\frac{\log 3^{-1}}{\log \theta}\big]\le j\le \frac{\log \delta^{-1/2}}{\log \theta}}
\exp\{-\frac{x^2}{2}\} \exp\left\{ \left(\frac{1}{2}\theta^{2j+2} - \frac{1}{2}\theta^{2j}\right)x^2+O_{\theta^{j+1},\theta^{2j}/2}x^2\Delta_{n,x}\right\}\\
\le & \sum_{\big[\frac{\log 3^{-1}}{\log \theta}\big]\le j\le \frac{\log \delta^{-1/2}}{\log \theta}}
 \exp\left\{-\frac{x^2}{2}+  \frac{x^2}{2}(\theta^2-1)\theta^{2j} +C\theta^{3j} e^{\theta/2}x^2\Delta_{n,x}\right\}\\
\le & \left(\frac{ \log \delta^{-1/2}  +\log 3}{\log \theta}+1\right)\exp\left\{-\frac{x^2}{2}+\frac{1}{2}(\theta^2-1)\delta^{-1}x^2+
C\delta^{-3/2} e^{\theta/2}x^2\Delta_{n,x}\right\}.
\end{align*}
Let $\theta^2=1+\delta x^{-2}$. It is easily seen that
$$ \frac{ \log \delta^{-1/2}  +\log 3}{\log \theta}+1 \le \exp\left\{O(1) (\log \delta^{-1}+\log x)\right\}. $$
It follows that
$$\Capc\left(S_n\ge x V_n, \delta \overline{B}_n^2\le V_n^2\le 9 \overline{B}_n^2\right) \le
\exp\left\{-\frac{x^2}{2} + O(1)
 \delta^{-3/2} \left(\log x+x^2\Delta_{n,x}\right)\right\}. $$
 The proof is completed. $\Box$

 Now, (\ref{eqth6.1.1}) and the upper bound of (\ref{eqth6.1.3}) follows from Propositions \ref{prop6.1}-\ref{prop6.3} immediately. As for the lower bound of (\ref{eqth6.1.3}), we let $b=b_n=x_n/\overline{B}_n$, $\eta_n=2b_n S_n-(b_n V_n)^2$.   Then by Lemma \ref{lem6.2} (a) and (c),
$$ \lim_{n\to \infty} x_n^{-2}\ln \Sbep\exp\left\{\lambda\eta_n\right\}=2\lambda^2-\lambda, \;\; \lambda>1/2, $$
which implies the following proposition similar to Proposition \ref{prop4}.
\begin{proposition}\label{prop6.4} For any $0<\beta<1$,
\begin{equation}\label{eqprop6.4.1}
\liminf_{n\to \infty} x_n^{-2} \ln \Capc\left( \eta_n\ge (1+\beta)x_n^2\right)\ge -\frac{(1+\beta/4)^2}{2}.
\end{equation}
\end{proposition}

Then, the lower bound of (\ref{eqth6.1.3}) follows by noting
$$ \Capc(S_n\ge x_nV_n)\ge    \Capc( 2b_n  S_n-(b_nV_n)^2\ge x_n^2). $$
The proofs are now completed.

\end{document}